\documentclass[12pt,a4paper,reqno]{amsart}
\usepackage{amsmath}
\usepackage{amsfonts}
\usepackage{amssymb}

\numberwithin{equation}{section}

\newtheorem*{main}{Main Theorem}
\newtheorem{theorem}{Theorem}[section]
\newtheorem{corollary}{Corollary}[section]
\newtheorem{lemma}{Lemma}[section]

\newtheorem{definition}{Definition}[section]

\def\A{{\mathbb A}}

\def\F{{\mathbb F}}

\def\fG{{\mathfrak G}}
\def\H{{\mathbb H}}

\def\K{{\mathfrak K}}
\def\L{{\mathfrak L}}

\def\N{{\mathcal N}}
\def\O{{\mathcal O}}

\def\Q{{\mathbb Q}}

\def\cS{{\mathfrak S}}
\def\W{{\mathfrak W}}
\def\cW{{\mathcal W}}
\def\T{{\mathcal T}}
\def\U{{\mathcal U}}

\def\Y{{\mathcal Y}}
\def\Z{{\mathbb Z}}

 \DeclareMathOperator{\Cl}{Cl}
\DeclareMathOperator{\coker}{coker} \DeclareMathOperator{\Div}{Div}
\DeclareMathOperator{\Gal}{Gal} 
 \DeclareMathOperator{\image}{Im}
 
\DeclareMathOperator{\Hom}{Hom} \DeclareMathOperator{\ord}{ord}
 \DeclareMathOperator{\Prin}{P}

\begin{document}
\title{On Iwasawa Theory over Function Fields}

\author[Kueh]{Ka-Lam Kueh$^+$}
\thanks{$^+$ Prof. Kueh passed away on July 24th 2005}
\address{Institute of mathematics\\
Academia Sinica\\
Taipei 11529, Taiwan}
\author[Lai]{King Fai Lai}
\address{School of Mathematics and Statistics\\
    University of Sydney\\
    NSW 2006, Australia}
\email{kflai@math.usyd.edu.au}
\author[Tan]{Ki-Seng Tan$^{\times}$}
\address{Department of Mathematics\\
National Taiwan University\\
Taipei 10764, Taiwan} \email{tan@math.ntu.edu.tw}

\thanks{$^{\times}$ The author was supported in part by the National
Science Council of Taiwan (R.O.C.), NSC91-2115-M-002-001,
NSC94-2115-M-002-010, NSC95-2115-M-002-017-MY2 }



\maketitle

\begin{abstract}
Let $k_{\infty}$ be a $\Z_p^d$-extension of
 a global function field $k$ of characteristic $p$.
Let $\Cl_{k_{\infty},p}$ be the $p$ completion of the class group of
$k_{\infty}$. We prove that the characteristic ideal of the
Galois module $\Cl_{k_{\infty},p}$ is generated by the Stickelberger element
of Gross which calculates the special values of $L$ functions
\end{abstract}

\begin{section}{Introduction}\label{s:intro}

We prove the $p$ part of  Iwasawa's main
conjecture over global function fields. We show that the
characteristic ideal of the Iwasawa module obtained from the class
groups is generated by the Stickelberger element defined by Gross.
Our proof begins with a theorem of Tate and the rest of the proof
is based on some results on the leading term of
the Stickelberger element first conjectured by Gross (\cite{g1}) and on
the properties of pseudo-null Iwasawa modules studied by Greenberg
(\cite{grn}).

Let $k$ be a global function field of
characteristic $p$. Let $k_{\infty}/k$ be a $\Z_p^d$-extension unramified outside a
non-empty finite set
$S$ of places of $k$.
For each non-negative integer $n$
let $k_n$ be the fixed field of $\Gamma^{p^n}$.
 Denote the $p$-Sylow subgroup of the class group $\Cl_{k_n}$ of
 divisors of degree zero
 of $k_n$ by $\Cl_{k_n,p}$.
Let $\Cl_{k_{\infty},p}$ be the projective limit
$\lim_{\stackrel{\longleftarrow}{n}} \Cl_{k_n,p}$
taken over
 the norm
maps  $\Cl_{k_m,p}\longrightarrow
\Cl_{k_n,p}$ with $m\geq n$.

Write $\Gamma$
  for the Galois group $\Gal(k_{\infty}/k)$.
Let $\Lambda_{\Gamma}$ denote the complete group ring $\Z_p[[\Gamma]]$.
Let $\sigma_1,...,\sigma_d$ be a
$\Z_p$-basis of $\Gamma$ viewed as a $\Z_p$-module.  Put
$t_j=\sigma_j-1$.
Then $\Lambda_{\Gamma}=\Z_p[[t_1,...,t_d]]$.
The ring
$\Lambda_{\Gamma}$ is a unique factorization domain (\cite{bou},
Chap. 7, Sec.3.9, Prop. 8) and
the class group $\Cl_{k_{\infty},p}$ is a finitely generated torsion
$\Lambda_{\Gamma}$ module (\cite{Iw59,ser}).
Moreover $\Cl_{k_{\infty},p}$ is pseudo-isomorphic to
$\oplus_{j=1}^J \Lambda_{\Gamma}/\wp_j^{n_j}$
where the  non-negative integers $n_j$ and  the height $1$
 prime ideals $\wp_j$ of
$\Lambda_{\Gamma}$ are uniquely determined
(\cite{bou}, Chap.
VII, Sec. 4.4, Definition 3 and Theorem 5).
We define
the characteristic ideal $\chi_{{\Gamma}}(\Cl_{k_{\infty},p})$
of the $\Lambda_{\Gamma}$-module by
\begin{equation}\label{e:charide}
\chi_{{\Gamma}}(\Cl_{k_{\infty},p}) :=
\prod_{j=1}^J \wp_j^{n_j}.
\end{equation}

Fix another finite set  $T$ of places of $k$ so that $S\cap T=\emptyset$.
Given a finite
abelian extension $K/k$ unramified outside $S$ with Galois group
$G$, let
$[v]\in G$ be the Frobenius at $v$ and $q$ be  the order of the
constant field of $k$.
 Put
\begin{equation}\label{e:thetau}
\Theta_{K/k,S,T}(u)=
\prod_{v\notin S}(1-[v]u^{\deg v})^{-1} \cdot \prod_{v\in T}(1-[v](q
u)^{\deg v}).
\end{equation}
Let $\theta_{K/k,S,T}=\Theta_{K/k,S,T}(1)$.
The result of Weil (see \cite{t1} Chapter V) shows that
$\Theta_{K/k,S,T}(u)$ is a polynomial in $u$ with coefficients in
$\Z[G]$, in particular $\theta_{K/k,S,T} \in \Z[G]$ and is the same as
the Stickelberger element $\theta_G$ defined in \cite{g1}.

If $K'$ is an intermediate field of $K/k$
with $\Gal(K/K')=H$, then the Stickelberger elements enjoy the
functorial property that under the ring homomorphism
$pr:\Z[G]\longrightarrow \Z[G/H]$ induced from the natural
projection from $G$ onto $G/H$ we have
\begin{equation}\label{e:functo}
pr(\theta_{K/k,S,T})=\theta_{K'/k,S,T}.
\end{equation}
We use these maps to define $\theta_{k_{\infty}/k,S,T}\in\Lambda_{\Gamma}$
as the projective limit over $n$ of $\theta_{k_n/k,S,T}\in\Z_p[\Gamma_n]$.

We can now state our main result.

\begin{main}\label{t:main}
Let $k$ be a global function field of characteristic $p$. Let
$k_{\infty}/k$ be
  a $\Z_p^d$-extension with Galois group
$\Gamma=\Gal(k_{\infty}/k)$. We assume that  $k_{\infty}/k$ is unramified
outside a non-empty finite set $S$ of places $k$ and that no place
in $S$ splits completely in $k_{\infty}/k$. Let $S_0\subset S$
denote the subset consisting of unramified places in $k_{\infty}/k$.
If $S_0\not=S$, put
$$\delta=\prod_{v\in S_0} (1-[v]),$$
where $[v]\in \Gamma$ is the Frobenius at $v$; if
$S_0=S$, then $k_{\infty}/k$ is the constant $\Z_p$-extension, and
we put
$$\delta=(1-Fr)^{-1}\prod_{v\in S_0} (1-[v]),$$
where $Fr\in \Gamma$ is the Frobenius element. Fix a
non-empty finite set $T$  of places outside $S$. Then
$$(\theta_{k_{\infty}/k,S,T})
=(\delta)\cdot \chi_{\Gamma}(\Cl_{k_{\infty},p}).$$
\end{main}

\noindent
{\bf Remarks}: \quad
Over a global function field of characteristic different from $p$,
the constant field extension is the only $\Z_p^d$-extension. In this
case $d=1$ and Iwasawa's main conjecture is known to be a
consequence of Weil's theory. Therefore, we only deal with
 those $\Z_p^d$-extensions with $p$ equals the characteristic of
the ground field. Also,we shall see in Section \ref{subs:rr} that,
unlike the counterpart in the cases where $k$ is a number field,
for a fixed $S$ the rank $d$ of the Galois group $\Gamma$ can be
arbitrary large. The next point is that the assumption in the
theorem is necessary, since we know that if some $v\in S$ splits
completely then $\theta_{k_{\infty}/k,S,T}=0$
(\cite{g1},\cite{tn}). The appearance of the extra factor
$\delta$ is  consistent with the fact that if $v\in S_0$ and
$S'=S-\{v\}$, then
\begin{equation}\label{e:compare}
\theta_{k_{\infty}/k,S,T}=(1-[v])\theta_{k_{\infty}/k,S',T}.
\end{equation}
We  note that if $S_0=S$, then
\begin{equation}\label{e:frv}
[v]=Fr^{\deg_k(v)}
\end{equation}
and hence
$(1-Fr)^{-1}(1-[v])\in\Lambda_{\Gamma}$.
Finally, if a different $T'$ is chosen, then the new
Stickelberger element $\theta_{k_{\infty}/k,S,T'}$
will be a multiple of $\theta_{k_{\infty}/k,S,T}$ by a unit in $\Lambda_{\Gamma}$,
and the two Stickelberger elements generate the same ideal. This explains the dependence on $T$ in
the last equation in the main Theorem.

Here is an outline of the paper. The proof of the Main theorem will
be completed in section 4. In section 3 we set up the induction machine which
allows us to bootstrap one $\Z_p$ extension down each step. This is summarized in
Lemma \ref{l:finally}.
In section 2
we prove the results which are needed to begin the induction, namely we use Tate's theorem to prove the
Main Theorem in the case of a $\Z_p$ extension whose the Stickelberger element is monomial
(in the sense of Definition \ref{d:mono})
and we construct an independent extensions (see Definition \ref{d:indp})
$L_{\infty}$ over $k_{\infty}$ which contains
an intermediate $\Z_p$ extension whose
Stickelberger element is monomial.
We shall keep the notations of this section in the rest of the paper.

\end{section}

\begin{section}{The class groups and the regulators}\label{sec:gen}

In this section we first  draw some consequences from Tate's theorem.
Next  we
recall Weil's theory on the zeta-function
associated to a global function field, and then we apply it together
with Tan's theorem \cite{tn} on
the exact form of the leading term
of the
Stickelberger element $\theta_{k_{\infty}/k,S,T}$
 as conjectured by Gross ( \cite{g1})
 and the computation of $p$-adic values of special values of $L$ function
 to prove the Main
Theorem  in the  case when the Stickelberger element
is monomial in the sense of Definition \ref{d:mono}.

We shall see that by enlarging the set $S$ of places of $k$
we can make the order of the class group $\Cl_{k,S,T}$ relatively prime to $p$.
It seems that we cannot simplify the refined regulator by just
enlarging the set $S$. But we can do so by extending the field
$k_{\infty}$.
The other main result of this section is the existence of an independent
extension of the given pair $(k,S)$ (see Definition \ref{d:indp} and
Lemma \ref{l:indp}). Then it follows immediately from Tate's theorem that
if $(L_{\infty}, {\tilde S})$ is an independent extension of
$(k_{\infty}, S )$ then
\begin{equation}\label{e:indepen}
\chi_{{\Gal(L_{\infty}/k)}}(Cl_{L_{\infty},p})
=(\theta_{L_{\infty}/k,{\tilde S},T}^m),\;
\end{equation}
for some\;positive integer $ m$.

The group of divisors of $k$, denoted by $\Div_k$,
is the free $\Z$-module
generated by all the places of $k$ ( see \cite{t1} Chap. V). There is a group homomorphism
$k^{\times}\longrightarrow \Div_k$ such that an element $\alpha\in k^{\times}$ is sent to the divisor
$\sum_v \ord_v(\alpha)\cdot v$, and the image of this homomorphism is denoted as $\Prin_k$.
Let $\Div_k^0$ denotes the subgroup of $\Div_k$ consisting of
zero-degree divisors. The class group of $k$ is
  $\Cl_k : = \Div_k^0/\Prin_k$.

\begin{subsection}{Some consequences of Tate's theorem}\label{subs:tate}

 Our proof of the Main Theorem
begins with a theorem of Tate (\cite{t1}), which is in fact the
function-field version of the Brumer-Stark conjecture.
We shall only quote the main part of this splendid theorem.
Let $S$ and $T$ be as before.
For a finite
abelian extension $K/k$ unramified outside $S$ with Galois group
$G=\Gal(K/k)$, put
\begin{equation}\label{e:Thetau}
\Theta_{K/k,S}(u)=\prod_{v\notin S}(1-[v]u^{\deg v})^{-1}\in
1+u\Q[G][[u]],
\end{equation}
where $[v]\in G$ denotes the Frobenius at $v$
and  $q_K$ denotes the order of the
constant field $\F_K$ of $K$.

\begin{theorem}\label{t:tate} \em{(Tate \cite{t1})}
We have $(q_K-1)\cdot \Theta_{K/k,S}(1)\in \Z[G]$ and
$$(q_K-1)\cdot \Theta_{K/k,S}(1)\cdot \Cl_K=0.$$
\end{theorem}

Obviously, we have
\begin{equation}\label{e:thetau}
\Theta_{K/k,S,T}(u)=\Theta_{K/k,S}(u)\cdot \prod_{v\in T}(1-[v](q
u)^{\deg v}).
\end{equation}
Since the number $q_K-1$ is relatively prime to $p$, it follows from Tate's
theorem and the equation (\ref{e:thetau}) that
\begin{equation}\label{e:ann}
\theta_{K/k,S,T}\cdot\Cl_{K,p}=0
\end{equation}
and so
\begin{equation}\label{e:anninfty}
\theta_{k_{\infty}/k,S,T}\cdot\Cl_{k_{\infty},p}=0.
\end{equation}
Since $\Lambda_{\Gamma}$ is a unique factorization domain, every prime
ideal $\wp_i$ appearing in the characteristic ideal
$\chi_{{\Gamma}}(\Cl_{k_{\infty},p})$ (\ref{e:charide})
 is generated by a prime element
$\pi_i\in\Lambda_{\Gamma}$. Therefore, the characteristic ideal
$\chi_{{\Gamma}}(\Cl_{k_{\infty},p})$ is generated by the
product
$\varpi:=\prod_{i\in J}\pi_i^{n_i}$.
Now the equation (\ref{e:anninfty}) is actually equivalent to saying that
$\theta_{k_{\infty}/k,S,T}$ is divisible by the least common
multiple of all the $\pi_i^{n_i}$.

\begin{corollary}\label{c:tate}
If $\theta_{k_{\infty}/k,S,T}$ is irreducible in $\Lambda_{\Gamma}$,
then
\begin{equation}\label{e:cor1.1}
\chi_{{\Gamma}}(\Cl_{k_{\infty},p})=(\theta_{k_{\infty}/k,S,T}^m),
\end{equation}
for some non-negative integer $m$.
\end{corollary}

We will prove that for suitable choice of $S$, the Stickelberger
element $\theta_{k_{\infty}/k,S,T}$ actually gives rise to a
generator of the characteristic ideal of $\Cl_{k_{\infty},p}$.

\end{subsection}

\begin{subsection}{The class number formula of Gross}\label{subs:zeta}
Consider a finite extension
$K/k$ unramified outside $S$.
For each place $w$ of $K$ let $\F_w$ denote the residue field of $w$.
Put $q_w=|\F_w|$. The zeta-function
$\zeta_{K}(s)=\prod_{w} (1-q_w^{-s})^{-1}$
can be written as
\begin{equation}\label{e:weil}
\zeta_K(s)=\frac{P_1(q_K^{-s})}{(1-q_K^{-s})(1-q_K^{1-s})}
\end{equation}
with $P_1(u)\in\Z[u]$ and
\begin{equation}\label{e:weilclass}
P_1(1)=|\Cl_K|.
\end{equation}

Use $S(K)$ (resp. $T(K)$) to denote the set consisting of places of
$K$ sitting over $S$ (resp. $T$), and define
\begin{equation}\label{e:thetastdef}
\zeta_{K,S,T}(s)=\prod_{w\in S(K)}{1-q_w^{-s}}\cdot\prod_{w\in T(K)}
(1-q_w^{1-s})\cdot\zeta_K(s).
\end{equation}
Recall that the degree $\deg_K(w)$ of $w$ equals the degree
$|\F_w:\F_K|$ of the field extension $\F_w/\F_K$. It is easy to see
from equalities (\ref{e:weil}) and (\ref{e:weilclass})
that if
$r_K=|S(K)|-1$,
then the Taylor expansion begins as
\begin{equation}\label{e:zeta}
\zeta_{K,S,T}(s)\equiv B_{K,S,T}\cdot s^{r_K} \pmod{s^{r_K+1}},
\end{equation}
where
\begin{equation}\label{e:bkst}
B_{K,S,T}=|\Cl_K| (\log q_K)^{r_K}(\prod_{w\in S(K)}
\deg_K(w))\prod_{w \in T(K)} (1-|\F_w|)
\end{equation}
For an intermediate field $K'$ of $K/k$ we define the ratio $\rho_{K/K'}(s)$
of zeta functions as
$$\rho_{K/K'}(s):=\zeta_{K,S,T}(s)/\zeta_{K',S,T}(s).$$
The next lemma is an immediate consequence of (\ref{e:zeta}) and (\ref{e:bkst}).

\begin{lemma}\label{l:ratio}
If $K'$ is an intermediate field of $K/k$ with $r_{K'}=r_K$, then
the function $\rho_{K/K'}(s)$ is
regular at $s=0$ with
$$\rho_{K/K'}(0)=\epsilon_0\cdot \frac{|\Cl_K|}{|\Cl_{K'}|}
\cdot (\frac{\log q_K}{\log q_{K'}})^{r_K}
\cdot\frac{\prod_{w\in S(K)}\deg_K(w)}
{\prod_{w\in S(K')}\deg_{K'}(w)},$$
for an $\epsilon_0\in\Z_p^{\times}$.
\end{lemma}

We would like to point out that (\ref{e:zeta}) is a class number formula since
the coefficient $B_{K,S,T}$ can be expressed as the product of some
class number and regulator. To see this, we consider $\O_{K,S}$, the
ring of $S$-integers of $K$. Let $\Cl_{K,S}$ denote the ideal class
group of the ring $\O_{K,S}$ and let $U_{K,S,T}$ denote the subgroup
of $\O_{K,S}^{\times}$ consisting of elements which are congruent to $1$
modulo $T(K)$. Define $\Div_{K,S,T}$ as the group whose elements are
of the form $( \mathcal{I}, (a_w)_{w\in T(K)})$ where $ \mathcal{I}$
is an ideal of $\O_{K,S}$ and $a_w$ is a local generator of $
\mathcal{I}$  at $w$. And define $\Cl_{K,S,T}$ as the quotient of
$\Div_{K,S,T}$ by the subgroup generated by principle elements. Then
we have the exact sequence (\cite{g1})
\begin{equation}\label{e:exact}
1\longrightarrow U_{K,S,T}\longrightarrow \O^{\times}_{K,S}\longrightarrow
\prod_{w\in T(K)}\F_w^{\times}\longrightarrow \Cl_{K,S,T}\longrightarrow
\Cl_{K,S}\longrightarrow 1.
\end{equation}
Let $R_{K,S,T}$ denote the classical regulator formed by
$U_{K,S,T}$. Namely, if $w_1,...,w_{r_K}$ are distinct places in
$S(K)$ and $u_1,...,u_{r_K}$ is a $\Z$-basis of $U_{K,S,T}$, which
is in fact free over $\Z$ (see \cite{g1}), then
$R_{K,S,T}=|{\det}_{1\leq i,j\leq r_K} (\ord_{w_i}(u_j)\cdot\deg_K(w_i))| $.
This implies (see \cite{g1} )
\begin{equation}\label{e:bkst2}
B_{K,S,T}=  (-1)^{|T(K)|-1}\cdot {|\Cl_{K,S,T}|\cdot R_{K,S,T}}.
\end{equation}

We conclude this section by recalling another type of class number
formula proposed by Gross. It involves the Stickelberger element and
a refined regulator (\cite{g1}). This formula will be useful for us.

Assume that $K/k$ is a pro-$p$ abelian extension with $G=\Gal(K/k)$. If
$G$ is finite, $I_G$ (resp. $I_{G,p}$), the augmentation ideal of
$\Z[G]$ (resp. $\Z_p[G]$), is defined as the kernel of the ring
homomorphism onto $\Z$ (resp. $\Z_p$) sending $\sum_{g\in G}
a_g\cdot g$ to $\sum_{g\in G} a_g$. If $G$ is pro-finite and $n\in\Z_+$, the
$n$'th power
augmentation ideal $I_G^n$ (resp. $I_{G,p}^n$) is defined as the
projective limit of the corresponding $n$'th powers $I_{\bar
G}^n$ (resp. $I_{{\bar G},p}^n$) where ${\bar G}$ runs
through all finite quotients of $G$.

For a place $v\in S$ let
\begin{equation}\label{e:localnorm}
\psi_{v,G}:k_v^{\times}\longrightarrow G_v\subset G
\end{equation}
be the local norm residue map.
Write $S=\{v_0,v_1,...,v_r\}$ and $r=r_k=|S|-1$. We choose a
$\Z$-basis $u_1,...,u_r$ of $U_{k,S,T}$. Assume that the ordering of
this basis is chosen such that the number
$${\det}_{1\leq i,j\leq r} (\ord_{v_i}(u_j)\cdot\deg_k(v_i))$$
is positive. Then the refined regulator of Gross is defined as the
residue class modulo $I_{G}^{r+1}$ of the element
\begin{equation}\label{e:greg}
{\det}_{K/k,S,T}=\det_{1\leq i, j \leq
r}(\psi_{v_i,G}(u_j))\in I_G^r\subset \Z[[G]].
\end{equation}

The following pro-$p$ version of a conjecture of Gross (\cite{g1})
was proved in \cite{tn} (the case where $r=1$ was first proved in
\cite{h}).

\begin{theorem}\label{t:gross}
If $K/k$ is a pro-$p$ abelian extension unramified outside $S$, then
$\theta_{K/k,S,T}\in I_{G,p}^{r}$ and
$\theta_{K/k,S,T}\equiv |\Cl_{K,S,T}|\cdot {\det}_{K/k,S,T}
\pmod{I_{G,p}^{r+1}}$.
\end{theorem}

\end{subsection}

\begin{subsection}{Special values of $L$ function}\label{subs:pv}

Let  $S$ and $T$ be as before.
Write $q_v=q^{\deg v}$ for the order of the residue field at a place $v$ of $k$.
Let $K/k$ be  a finite
abelian extension  unramified outside $S$ with Galois group
$G$. Write
$[v]\in G$ for the Frobenius at $v$.
For each character $\varphi$ of $G$  define
$$L_{\varphi,S,T}(s)=\prod_{v\not\in S}\frac{1}{1-\varphi([v])q_v^{-s}}\cdot\prod_{v\in T}
(1-\varphi([v])q_v^{1-s}).$$
We extend $\varphi$  linearly to a
homomorphism $\varphi:\Z[G][[q^{-s}]]\longrightarrow {\overline
{\Q}}[[q^{-s}]]$. It follows from (\ref{e:thetau}) that
$$\varphi(\Theta_{K/k,S,T}(q^{-s}))= L_{\varphi,S,T}(s).$$

If we view the dual group ${\widehat {G/H}}$ as a subgroup of $\widehat G$ then the ratio
$\rho_{K/K'}(s)$ of zeta functions can be written as a product of $L$ functions, namely,
$
\rho_{K/K'}(s)=\prod_{\varphi\not\in {\widehat {G/H}}} L_{\varphi,S,T}(s),
$
and so
\begin{equation}\label{e:prodtheta}
\prod_{\varphi\not\in {\widehat
{G/H}}}\varphi(\Theta_{K/k,S,T}(q^{-s}))=\rho_{K/K'}(s).
\end{equation}

Now Lemma \ref{l:ratio} implies the following.

\begin{lemma}\label{e:thetaratio}
If $K'$ is an intermediate field of $K/k$ with $r_{K'}=r_K$, then
$$\prod_{\varphi\not\in {\widehat {G/H}}}
\varphi(\theta_{K/k,S,T})=\epsilon_0\cdot \frac{|\Cl_K|}{|\Cl_{K'}|}
\cdot (\frac{\log q_K}{\log q_{K'}})^{r_K} \cdot\frac{\prod_{w\in
S(K)}\deg_K(w)} {\prod_{w\in
S(K')}\deg_{K'}(w)},\;\;\epsilon_0\in\Z_p^{\times}.$$
\end{lemma}

A continuous character $\varphi\in {\widehat {\Gamma}}$ can also be
extended linearly to a ring homomorphism
$\varphi:\Lambda_{\Gamma}\longrightarrow {\overline \Q}_p$. If
$\varphi\in{\widehat \Gamma}_n \subset {\widehat {\Gamma}_{\infty}}$
then it has order dividing $p^n$, and vice versa. From the
functorial property (\ref{e:functo}), we have
$\varphi(\theta_{k_{\infty}/k,S,T})=\varphi(\theta_{k_n/k,S,T})$.

Now we assume that $k_{\infty}/k$ is a $\Z_p$-extension and
let $\ord_p$ be the unique valuation on ${\bar \Q}_p$ with
$\ord_p(p)=1$. If $\varphi,\varphi' \in {\widehat
{\Gamma}_{\infty}}$ are of the same order, then they are conjugate
under the action of $\Gal({\bar \Q}/\Q)$ and hence
\begin{equation}\label{e:v}
\ord_p(\varphi'(\theta_{k_{\infty}/k,S,T}))
=\ord_p(\varphi(\theta_{k_{\infty}/k,S,T})).
\end{equation}

Let
 $\sigma$ be a
$\Z_p$-basis of the Galois group
$\Gamma$ of $k_{\infty}/k$ viewed as a $\Z_p$-module.  Put
$t=\sigma-1$.
Then $\Lambda_{\Gamma}=\Z_p[[t]]$.
The prime element
$\pi_i$ generating the
 prime
ideal $\wp_i$ appearing in the characteristic ideal
$\chi_{{\Gamma}}(\Cl_{k_{\infty},p})$ (\ref{e:charide})
can be chosen to be
either the prime number $p$ or an irreducible
 distinguished polynomial $p_i(t)\in\Lambda_{\Gamma}$.
In this case we have the following  pseudo-isomorphism
\begin{equation}\label{e:elementary}
\Cl_{k_{\infty},p}\longrightarrow
\bigoplus_{i=1}^M\Lambda_{\Gamma}/p^{m_i}\Lambda_{\Gamma}\oplus
\bigoplus_{j=1}^L\Lambda_{\Gamma}/p_j(t)^{l_j} \Lambda_{\Gamma}.
\end{equation}
If we put
\begin{equation}\label{e:mulambda}
\mu=\sum_{i=1}^M m_i,\;\;f=\prod_{j=1}^L p_j(t)^{l_j}\;\;
\text{and}\;\; \lambda=\sum_{j=1}^L l_j\cdot\deg(p_j(t)),
\end{equation}
then (\cite{ser,Iw59})
there is a constant $\nu$ such that $|\Cl_{k_n,p}|=p^{e_n}$ with
\begin{equation}\label{e:iwasawa}
e_n=\mu p^n+\lambda n+\nu.
\end{equation}

\begin{lemma}\label{l:mulambda}
Assume that $k_{\infty}/k$ is a $\Z_p$-extension and no place in $S$
splits completely under $k_{\infty}/k$. Let $\mu$, $\lambda$, $f$ be
as in {\em (\ref{e:elementary})} and {\em (\ref{e:mulambda})} and
let $\delta$ be as in the Main Theorem. Then
$\theta_{k_{\infty}/k,S,T}=\epsilon_1\cdot p^{\mu}\cdot g$,
where $\epsilon_1\in\Lambda_{\Gamma}^{\times}$ and $g$ is a distinguished
polynomial of degree
$\lambda+\deg(\delta)$ which is equal to $\deg(f\cdot \delta)$.
\end{lemma}
\begin{proof}
If an $\eta\in\Lambda_{\Gamma}$ is expressed as
$$\eta=\epsilon_2\cdot p^{m}\cdot h,$$
where $\epsilon_2\in\Lambda_{\Gamma}^{\times}$ and
$h=t^l+b_{l-1}t^{l-1}+\dots + b_0$ is the associated distinguished
polynomial with $b_0,...,b_{l-1}\in p\cdot \Lambda_{\Gamma}$, then
for a $\varphi\in {\widehat {\Gamma}_{\infty}}$ of order $p^n$ with
$p^n\geq l$, we have
$$\ord_p(\varphi(\eta))=m+\frac{l}{p^n-p^{n-1}}.$$
In view of this, we only need to show that for $n$ large
$$\ord_p(\varphi(\theta_{k_{\infty}/k,S,T}))=\mu+\frac{\lambda+\deg(\delta)}{p^n-p^{n-1}}.$$ From
the equation (\ref{e:iwasawa}), Lemma \ref{e:thetaratio} and the
equation (\ref{e:v}), we see that it is enough to verify that for
large $n$
$$p^{\deg(\delta)}=(\frac{q_{k_n}}{q_{k_{n-1}}})^{r_{k_n}}\cdot \frac{\prod_{w\in S(k_n)}\deg_{k_n}(w)}
{\prod_{w\in S(k_{n-1})}\deg_{k_{n-1}}(w)}.$$ Since no place in $S$
splits completely, the cardinality of the set $S(k_n)$ will
eventually be stable, and every $w$ in $S(k_{n-1})$ will be either
ramified or inert under $k_{n}/k_{n-1}$. In the case where
$k_{\infty}/k$ is the constant $\Z_p$-extension and $n$ is large
enough, we have $q_{k_n}=q_{k_{n-1}}^p$, $r_{k_n}=\deg(\delta)$ and
$\deg_{k_n}(w)=\deg_{k_{n-1}}(w')$ if $w'$ sits below $w$. Therefore
the lemma holds. In other cases, when $n$ is large enough
$q_{k_n}=q_{k_{n-1}}$, for $w'$ sitting below $w$ the ratio
$\deg_{k_n}(w)/\deg_{k_{n-1}}(w')$ equals $1$ (resp. $p$) if
$k_n/k_{n-1}$ ramified (resp. inert) at $w'$, and the cardinality of
inert places in $S(k_{n-1})$ equals $\deg(\delta)$. The lemma also
holds in these situations.
\end{proof}

\begin{definition}\label{d:mono}
Assume that $k_{\infty}/k$ is a $\Z_p$-extension and no place in $S$
splits completely under $k_{\infty}/k$. The Stickelberger element
$\theta_{k_{\infty},S,T}$ is called monomial if for some
$\epsilon_2\in \Z_p^{\times}$
$$\theta_{k_{\infty},S,T}\equiv \epsilon_2 t^r\pmod{(t^{r+1})}.$$
\end{definition}

\begin{lemma}\label{c:simple}
Assume that $k_{\infty}/k$ is a $\Z_p$-extension and no place
in $S$ splits completely under $k_{\infty}/k$.
\begin{enumerate}
\item If the order of $\Cl_{k,S,T}$ is prime to $p$ and there exists a unit $\epsilon_2\in \Z_p^{\times}$
such that
$${\det}_{k_{\infty}/k,S,T}\equiv \epsilon_2\cdot t^r\pmod{(t^{r+1})},$$
then $\theta_{k_{\infty},S,T}$ is monomial.
\item If $\theta_{k_{\infty},S,T}$ is monomial, then
the Main Theorem  holds.
\end{enumerate}
\end{lemma}
\begin{proof}
The augmentation ideal $I_{\Gamma,p}$ is just the ideal generated by
$t$. The statement (1) is a consequence of Theorem \ref{t:gross}. If
the Stickelberger element is monomial, then
$\theta_{k_{\infty}/k,S,T}=t^r\cdot \xi$,
where the formal series $\xi$ begins with the constant term which is contained
in $\Z_p^{\times}$. This means that $\xi$ is itself a unit in
$\Lambda_{\Gamma}$. Lemma \ref{l:mulambda} says that $\mu=0$ and
$\deg(f)+\deg(\delta)=r$. From equations (\ref{e:anninfty}),
(\ref{e:compare}) and (\ref{e:frv}), we see that every prime factor
of $\varpi\delta$ divides $\theta_{k_{\infty}/k,S,T}$. Thus $t$ is
the only prime factor of both side.
\end{proof}

\end{subsection}

\begin{subsection}{Order of the group $\Cl_{k,S,T}$}\label{subs:clkst}

\begin{lemma}\label{l:clkst}
Let $\tilde S$ be a finite set of places of $k$ satisfying the
following conditions:
\begin{enumerate}
\item The subgroup of $\Cl_k$
generated by the set of all the zero-degree divisors which are
supported on $\tilde S$ contains the $p$-Sylow subgroup $\Cl_{k,p}$
of the class group.
\item The greatest common divisor of the degrees of the places in $\tilde S$ is one.
\end{enumerate}
Then the order of the group $\Cl_{k,{\tilde S},T}$ is prime to $p$.
\end{lemma}
\begin{proof}
Let $\mathbf{X}_{k,{\tilde S}}\subset\Div^0_k$ be the subgroup
formed by divisors supported on $\tilde S$. Use $\Div_{k,{\tilde
S}}$ to denote the group of ideals of the ring $\O_{\tilde S}$. Then
we have the exact sequence
\begin{equation}\label{e:xdiv}
0\longrightarrow \mathbf{X}_{k,{\tilde
S}}\stackrel{i}{\longrightarrow} \Div^0_{k}
\stackrel{\pi}{\longrightarrow}
\Div_{k,{\tilde S}}
\longrightarrow 0
\end{equation}
with $\pi$ taking $\sum_{v}a_v\cdot v $ to $ \sum_{v\not\in {\tilde S}} a_v\cdot v$.
The surjectivity of $\pi$ is from the condition (2). Consequently,
we have the induced exact sequence
\begin{equation}\label{e:xcl}
0\longrightarrow {\overline {\mathbf{X}}}_{k,{\tilde
S}}\stackrel{\bar i}{\longrightarrow} \Cl_{k}
\stackrel{\pi}{\longrightarrow}
\Cl_{k,{\tilde S}} \longrightarrow 0,
\end{equation}
where ${\overline {\mathbf{X}}}_{k,{\tilde S}}$ is the subgroup of
$\Cl_k$ formed by divisor classes obtained from $
\mathbf{X}_{k,{\tilde S}}$. The condition (1) says that $\bar i$ is
actually surjective on the $p$-part. Therefore $\Cl_{k,{\tilde S}}$
has order prime to $p$ and hence so is $\Cl_{k,{\tilde S},T}$ (see
(\ref{e:exact}) ).
\end{proof}

\begin{corollary}\label{c:clks}
There are infinitely many finite sets ${\tilde S}$ of
places of $k$ with the following properties:
\begin{enumerate}
\item the order of the group $\Cl_{k,{\tilde S},T}$ is prime to $p$,
\item the intersection ${\tilde S}\cap T=\emptyset$,
\item ${\tilde S}\supsetneqq S$ and no place in $\tilde S$ splits completely over $k_{\infty}/k$.
\end{enumerate}
\end{corollary}

\begin{proof}
Suppose the greatest common divisor of the degrees of the places in
$S$ is $N$. Let $L/k$ be the constant field extension of degree $N$.
Choose a generator $\tau\in\Gal(L/k)$.
Tchebotarev's density theorem says that there is a place $v$ of $k$
outside $S\cup T$ such that the element $\tau$ equals the Frobenius
$[v]\in\Gal(L/k)$.
We know from the class field theory that if $Fr_q:x\mapsto x^q$ is
the Frobenius substitution on the constant fields, then $Fr_q$
generate $\Gal(L/k)$ and $[v]=Fr_q^{\deg(v)}$. Since
$Fr_q^{\deg(v)}$ is also a generator of $\Gal(L/k)$, the degree of
$v$ must be relatively prime to $N$. Therefore the greatest common
divisor of the degrees of places in $S':=S\cup\{v\}$ equals $1$.

Consider the finite extension $L'/k$ which is the composite of all
the everywhere  unramified cyclic extensions of order $p$ and choose
a set of generators $\tau'_1,...,\tau'_s\in\Gal(L'/k)$.
Again, Tchebotarev's density theorem says that there are places
$v'_1,...,v'_s$ outside $S'\cup T$ such that each $\tau'_i$ equals
the Frobenius $[v'_i]\in\Gal(L'/k)$.
And we know from the class field theory that the Galois group
$\Gal(L'/k)$ is identified with $(\Div_k/\Prin_k)\otimes_{\Z}
\Z/p\Z$ and the condition on $\tau'_1,...,\tau'_s$ implies that the
divisor classes of $v'_1,...,v'_s$ generate
$(\Div_k/\Prin_k)\otimes_{\Z} \Z/p\Z$. Therefore, the divisor
classes of $v'_1,...,v'_s$ generate the $\Z_p$-module
$(\Div_k/\Prin_k)\otimes_{\Z} \Z_p$. Take ${\tilde
S}=S'\cup\{v'_1,...,v'_s\}$. Then the classes of all zero-degree
divisors supported on $\tilde S$ generate $\Cl_{k}\otimes_{\Z}
\Z_p$. Therefore the conditions (1) and (2) of Lemma \ref{l:clkst}
are satisfied. Tchebotarev's theorem also ensures us that the places
$v,v_1'...,v_s'$ can be chosen so that none of them splits completely over $k_{\infty}/k$.
\end{proof}
\end{subsection}

\begin{subsection}{Local norm residue maps}\label{subs:rr}

Let $\tilde S\supseteq S$ be a finite set of places and let $\H
$ be the Galois group of the maximal pro-$p$ abelian extension over
$k$ unramified outside $\tilde S$. It is known that (\cite{k,tn})
$\H$ is
isomorphic, as a topological group, to a countable infinite product
of $\Z_p$. This is actually due to the following simple Lemma
which can be viewed as the function-field version of the local Leopoldt conjecture.

\begin{lemma}\label{l:llc}
If at some place $v$ a global element $a\in k^{\times}$ equals $b^p$ for
some $b\in k_v^{\times}$, then $b\in k^{\times}$.
\end{lemma}
\begin{proof}
Since the field extension $k(a^{\frac{1}{p}})/k=k(a^{\frac{1}{p}})\cap k_v/k$ is both
purely inseparable and separable.
\end{proof}
Suppose that $K/k$ is a $\Z_p$-extension unramified outside $\tilde
S$. Then choosing a topological generator of $\Gal(K/k)$ is the same
as choosing an isomorphism $\Gal(K/k)\simeq\Z_p$. Thus the extension
$K/k$ together with a topological generator of $\Gal(K/k)$ gives
rise to a continuous homomorphism $\varphi:\H\longrightarrow \Z_p$,
and vice versa. In particular, taking $\tilde S=\{v\}$, we see that
there exists a $\Z_p$-extension, associated to a
$\varphi(v):\H\longrightarrow \Z_p,$
which is ramified at $v$ and unramified at other places.
\begin{lemma}\label{l:varphi}
There exists a $\Z_p$-extension over $k$, which is unramified
outside $\tilde S$ but ramified at every place in $\tilde S$.
\end{lemma}
\begin{proof}
The $\Z_p$-extension associated to $\sum_{v\in\tilde S}\varphi(v)$
satisfies the required condition.
\end{proof}

Let ${\tilde r}+1$ be the cardinality of ${\tilde S}$.
Choose distinct places $v_1,...,v_{\tilde r}\in {\tilde S}$. As
in (\ref{e:localnorm}), for each $i$ let
$$\psi_{v_i,\H}:k_{v_i}^{\times}\longrightarrow \H_{v_i}\subset \H$$
be the local norm residue map. Let $\varphi_i$ be the composition
of the natural embedding $U_{k,{\tilde S},T}\to  k_{v_i}^{\times}$ with
${\psi_{v_i,\H}}$.
Since $\H$ is a
$\Z_p$-module and $\psi_{v,i}$ is continuous, we can extend linearly
$\varphi_i$ to a homomorphism
$\varphi_i:\U\longrightarrow \H$
where
$\U:=U_{k,{\tilde S},T}\otimes\Z_p.$

\begin{lemma}\label{l:directsum}
Let $\cW:=\U\times \U\times\dots\times \U$ be the direct sum of
$\tilde r$ copies of $\U$. Then the homomorphism $\Psi:\cW
\longrightarrow \H$ sending $(w_1,...,w_{\tilde r})\in\cW$ to
$\sum_{i=1}^{\tilde r} \varphi_{i}(w_i)\in\H$ is injective and the
quotient group $\H/\Psi(\cW)$ is torsion free.
\end{lemma}

\begin{proof}
We need to show that if $w_1,...,w_{\tilde r}$ are elements in
$U_{k,{\tilde S},T}$ with $\sum_{i=1}^{\tilde r} \varphi_{i}(w_i)$
divisible by $p$ in $\H$, then every $w_i$ is divisible by $p$ in
$U_{k,{\tilde S},T}$. Let $x=(x_v)_v\in \A_k^{\times}$ be the idele such
that $x_v=w_i$ if $v=v_i$; $x_v=1$ otherwise. If $\sum_{i=1}^{\tilde
r} \varphi_{i}(w_i)$ is divisible by $p$, then there are
$y\in\A_k^{\times}$, $\alpha\in k^{\times}$ and $z\in \prod_{v\not\in\tilde
S}\O_v^{\times}$ such that
\begin{equation}\label{e:xalpha}
x=\alpha\cdot y^p\cdot z.
\end{equation}
As ${\tilde S}$ contains ${\tilde r}+1$ elements,
there is a place $v_0\in {\tilde S}\setminus \{v_1,...,v_{\tilde r}\}$. Then we see from the equality
(\ref{e:xalpha}) that $\alpha$ is divisible by $p$ in $k_{v_0}^{\times}$,
and hence by Lemma \ref{l:llc} there is an element $\beta\in k^{\times}$ such that
$\alpha=\beta^p$. The equality (\ref{e:xalpha}) implies that each
$w_i$ is divisible by $p$ in $k_{v_i}^{\times}$ and hence also divisible by
$p$ in $U_{k,{\tilde S},T}$.
\end{proof}

The following is a consequence of Lemma \ref{l:directsum} and the
fact that $\H$ is the direct product of countable infinite many
copies of $\Z_p$.

\begin{corollary}\label{c:independent}
Let  ${\tilde S}$ be a set of places of $k$ with ${\tilde r}+1$ elements.
If $K/k$ is any given abelian extension unramified outside $\tilde
S$ with Galois group isomorphic to $\Z_p^{d_0}$ for some
non-negative integer $d_0$, then there exists a field extension
$L/K$ with the following properties:
\begin{enumerate}
\item
$L/k$ is also an abelian extension unramified outside $\tilde S$
with Galois group isomorphic to $\Z_p^{c}$ for some non-negative
integer $c$.
\item If $u_1,...,u_{\tilde r}$ is a $\Z$-basis of
$ U_{k,{\tilde S},T}$ and $\psi_{v_i,\Gal(L/k)}$ is the local norm
residue map at $v_i$, then the subset
$\{\psi_{v_j,\Gal(L/k)}(u_i)\; | \; 1\leq i,j\leq {\tilde r}\}$ of
$\Gal(L/k)$ is linearly independent over $\Z_p$ and it generates a
direct summand of $\Gal(L/k)$.
\end{enumerate}
\end{corollary}

\end{subsection}
\begin{subsection}{Independent extensions}\label{l:indep}
In this section, we show that if the set $S$ and the field extension
$k_{\infty}$ are enlarged in a suitable way then the Stickelberger
element will become irreducible.

\begin{definition}\label{d:indp}
Let $(K,S)$ be a pair where $K/k$ is a $\Z_p^{d_0}$-extension
unramified outside $S$. A pair $(L,{\tilde S})$ is said to be  an independent
extension of $(K,S)$ if the following conditions hold:
\begin{enumerate}
\item  $\tilde S$ is  a finite set of places of $k$,
satisfying  {\em Corollary \ref{c:clks}}.
\item The field extension $L/k$
is a $\Z_p^c$-extension which ramifies at every place in $\tilde S$ and satisfies
{\em Corollary \ref{c:independent}}.
\end{enumerate}
\end{definition}

\begin{lemma}\label{l:indp}
\begin{enumerate}
\item There exist independent extensions $(L_{\infty},{\tilde S})$ of
$(k_{\infty},S)$  with
 $\tilde S$ arbitrarily large.
\item If $(L_{\infty},{\tilde S})$ is an independent extension of $(k_{\infty},S)$,
then $\theta_{L_{\infty}/k,{\tilde S},T}$ is irreducible in
$\Lambda_{\Gal(L_{\infty}/k)}$.
\item If $(L_{\infty},{\tilde S})$ is an
independent extension of $(k_{\infty},S)$ then there is an
intermediate $\Z_p$-extension $k'_{\infty}/k$ of $L_{\infty}/k$ such that
$k'_{\infty}/k$  ramifies at every place in $\tilde S$ and the
Stickelberger element $\theta_{k'_{\infty},{\tilde S},T}$ is
monomial.
\end{enumerate}
\end{lemma}

\begin{proof}
The existence of arbitrarily large $\tilde S$ follows  from Corollary
\ref{c:clks}. Let $K'/k$ be a $\Z_p$-extension unramified outside
$\tilde S$ but ramified at every place of $\tilde S$ (see Lemma
\ref{l:varphi}). Replace $K$ by the composite $K'K$ if necessary,
and we can assume that $K/k$ is ramified at every place of $\tilde
S$. Then the existence of $L_{\infty}$ is from Corollary
\ref{c:independent}. This proves (1).

To prove (2) we first recall the notations in Corollary
\ref{c:independent}. For $1\leq i,j\leq {\tilde r}$, we set
$\sigma_{i,j}=\psi_{v_j,\Gal(L_{\infty}/k)}(u_i)$ and
$t_{i,j}=\sigma_{i,j}-1$. Corollary \ref{c:independent} (2) says the
set $\{\sigma_{i,j}\}_{1\leq i,j \leq {\tilde r}}$ can be extended
to a basis $\{\sigma_1,...,\sigma_{c}\}$ of $\Gal(L_{\infty}/k)$
over $\Z_p$. If $t_i=\sigma_i-1$, then the augmentation quotient
$I_{\Gal(L_{\infty}/k),p}^{\tilde
r}/I_{\Gal(L_{\infty}/k),p}^{{\tilde r}+1}$ can be identified with
the $\Z_p$-module of $\tilde r$-degree homogeneous polynomials in
$t_1,...,t_c$. The refined regulator ${\det}_{L_{\infty}/k,{\tilde
S}.T}$ determines a residue class in the above augmentation
quotient, and from (\ref{e:greg}) we see that this residue class is
identified as the polynomial ${\det}(t_{i,j})_{1\leq i,j\leq {\tilde
r}}$. It is well-known (see \cite{v}) that this polynomial is
irreducible. Corollary \ref{c:clks} says that the order of
$\Cl_{k,{\tilde S},T}$ is prime to $p$, and then Theorem
\ref{t:gross} says that the Taylor expansion of
$\theta_{L_{\infty},{\tilde S},T}\in\Z_p[[t_1,...,t_c]]$ begins
with the irreducible polynomial
$|\Cl_{k,{\tilde S},T}|\cdot {\det}(t_{i,j})$
in $\Z_p[t_1,...,t_c]$.
Suppose $\theta_{L_{\infty},{\tilde S},T}=\theta_1\theta_2$ in
$\Z_p[[t_1,...,t_c]]$ and the Taylor expansions of $\theta_1$ and
$\theta_2$ begin with the leading homogeneous polynomials
$\vartheta_1, \vartheta_2\in\Z_p[t_1,...,t_c]$. Then the product
$\vartheta_1\vartheta_2$ is irreducible in $\Z_p[t_1,...,t_c]$.
Therefore one of them, say $\vartheta_1$ must be in the units group
$\Z_p^{\times}$ of the polynomial ring $\Z_p[t_1,...,t_c]$. This implies
that $\theta_1$, beginning with a unit in its Taylor expansion, must
be a unit in $\Z_p[[t_1,...,t_c]]$. This argument shows that
$\theta_{L_{\infty},{\tilde S},T}$ is an irreducible element in
$\Lambda_{\Gal(L_{\infty}/k)}$.

The $\Z_p$-module $\Hom_{\Z_p}(\Gal(L_{\infty}/k),\Z_p)$ is
isomorphic to $\Z_p^c$. Every element in it is a continuous map with
respect to the pro-finite topologies on $\Gal(L_{\infty}/k)$ and
$\Z_p$. Also, the pro-finite topology on $\Z_p^c$ coincides with the
compact-open topology on $\Hom_{\Z_p}(\Gal(L_{\infty}/k),\Z_p)$.
Since the subset $\Z_p^{\times}$ is open in $\Z_p$, the subset
$\mathrm{O}$ of the group $ \Hom_{\Z_p}(\Gal(L_{\infty}/k),\Z_p)$ consisting
of those $\psi$ such that
$${\det}_{\psi}:={\det}(\psi(\sigma_{i,j}))_{1\leq i,j,\leq {\tilde r}}\in \Z_p^{\times}$$
is open in $\Hom_{\Z_p}(\Gal(L_{\infty}/k),\Z_p)$. Since
$\{\sigma_{ij}\}_{1\leq i,j\leq\tilde r}$ is a subset of the basis
$\{\sigma_1,...,\sigma_c\}$, there is at least one $\psi$ such that
$\psi(\sigma_{ii})=1$, $i=1,...,\tilde r$, and $\psi(\sigma_{ij})=0$
for $i\not=j$. And for this $\psi$ the determinant
${\det}_{\psi}=1$. Therefore, the open set $
\mathrm{O}\not=\emptyset$. At each place $v\in\tilde S$ let $
\mathrm{O}(v)\subset\Hom_{\Z_p}(\Gal(L_{\infty}/k),\Z_p)$ be the
subset consisting of those $\psi$ whose restriction to the inertia
subgroup at $v$ is non-zero. Since this inertia subgroup is
non-trivial, the set $ \mathrm{O}(v)$ is a nonempty open subset of
$\Hom_{\Z_p}(\Gal(L_{\infty}/k),\Z_p)$.

Let $\psi\in \cap_{v\in\tilde S} \mathrm{O}(v)\cap \mathrm{O}$ and
let $k'_{\infty}$ be the fixed field of the kernel of $\psi$. Then
$k'_{\infty}/k$ is ramified at each place in $\tilde S$ and
$\theta_{k'_{\infty}/k,{\tilde S},T}$ is monomial (Lemma \ref{c:simple}).
\end{proof}

\end{subsection}
\end{section}

\begin{section}{The maximal pseudo-null sub-module}\label{sec:pseudonull}

Suppose that $\K_{\infty}/k$ is a
$\Z_p^e$-extension with $\Gal(\K_{\infty}/k)=\Upsilon$ and $\Xi$  is a
 rank one $\Z_p$-submodule of $\Upsilon$ with $\Upsilon/\Xi\simeq
\Z_p^{e-1}$. Let $\K_{\Xi}$ be the fixed field of $\Xi$ and let
$pr:\Lambda_{\Upsilon}\longrightarrow\Lambda_{\Upsilon/\Xi}$ be as usual the
natural projection. We shall find a  relation between
$pr(\chi_{{\Upsilon}}(\Cl_{\K_{\infty},p}))$ and
$\chi_{{\Upsilon/\Xi}}(\Cl_{\K_{\Xi},p})$ by using the fact that the
characteristic ideals are multiplicative (see \cite{bou}, Chap. 7,
Sec. 4.5).

Fix an
element $s\in\Lambda_{\Xi}$ such that $\xi:=s+1$ is a topological
generator of $\Xi$.
For an abelian Galois group $G$,
we use $G_v$, $G_{v}^0$ to denote the decomposition subgroup and the
inertia subgroup at a place $v$.
We choose a set $\cS$ of places of $k$ so that
$\K_{\infty}/k$ is unramified outside of $\cS$ and
no place in $\cS$ splits
completely in $\K_{\infty}/k$.
Let $\cS_1=\{v\in \cS \;\mid\;\Upsilon_v^0 \cap\Xi\not=\{0\}\}
$. Then $\K_{\infty}/\K_{\Xi}$ is unramified outside $\cS_1$.
We set
$\cS_2=\{v\in \cS_1\;\mid\; |   \Upsilon_v/\Upsilon_v^0|<\infty\}$.
Thus, a place $v\in\cS$ is in $\cS_2$ if and only if the extension
$\K_{\infty}/\K_{\Xi}$ is ramified at every place sitting over $v$ and the
corresponding residue field extension for $\K_{\infty}/k$ is of
finite degree.

\begin{subsection}{The  group  $\W_{\K_{\infty},p}$ }\label{subsec:clgl}


For a global function field $K$ let $\A_K^{\times}$ be  the idele group of
$K$.  Write $\W_K$ for the group $K^{\times}\backslash\A_K^{\times}/\prod_v\O_{K,v}^{\times}$.
Let
$\W_{K,p}$ be the $p$-completion of $\W_K$.
We define $\W_{\K_{\infty},p}$ as the projective limit
\begin{equation}\label{e:wkinfp}
\W_{\K_{\infty},p}=\lim_{\stackrel{\leftarrow}{K}} \W_{K,p},
\end{equation}
where
$K$ runs through all finite intermediate extensions
of $\K_{\infty}/k$.

The group  $\W_K$
is the group of divisor classes $\Div_K/\Prin_K$.
We have the
exact sequence
$0\longrightarrow \Cl_K\longrightarrow \W_K \stackrel
{\deg_K}{\longrightarrow} \Z\longrightarrow 0$,
where $\deg_K$
is the degree map. From the class field theory we see that $\W_K$ is
a dense subgroup of the Galois group of the maximal unramified
abelian extension of $K$ and the fixed field of $\Cl_K$ is the
maximal constant field extension of $K$. Therefore the $p$-completion
$\W_{K,p}$ of $\W_K$ is the Galois group of the maximal unramified
pro-$p$ abelian extension of $K$ and the fixed field of $\Cl_{K,p}$
is the constant $\Z_p$-extension of $K$. Since we also have
$\W_{K,p}=\W_K\otimes_{\Z}\Z_p$, an element in $\W_{K,p}$ can be
viewed as an equivalent class of $\Z_p$-divisors, which are $\Z_p$
linear combinations of places in $K$, where two such divisors are
equivalent if and only if they are $\Z_p$-linear equivalent in the
sense that they differ by a divisor of some element in the
$p$-completion ${\widetilde {K^{\times}}}$
of the multiplicative group $K^{\times}$. We should
note that if ${\widetilde { {\Prin}_K}}$ is the $p$-completion of the group
$\Prin_K$, then from the exact sequence
$$0\longrightarrow \F_K^{\times}\longrightarrow K^{\times}\longrightarrow \Prin_K\longrightarrow 0$$
and the fact that $|\F_K^{\times}|$ is prime to $p$, we get an isomorphism
${\widetilde {K^{\times}}} \simeq  {\widetilde { {\Prin}_K}}$.

If $K'$ is another finite intermediate field
containing $K$, then we have the commutative diagram:
$$
\begin{array}{ccccccc}
0\longrightarrow & \Cl_{K',p} & \longrightarrow & \W_{K',p} &
\stackrel{\deg_{K'}}{\longrightarrow} & \Z_p & \longrightarrow 0\\
{} & \downarrow & {} & \downarrow & {} & \downarrow & {}\\
0\longrightarrow & \Cl_{K,p} & \longrightarrow & \W_{K,p} &
\stackrel{\deg_{K}}{\longrightarrow} & \Z_p & \longrightarrow 0\\
\end{array}
$$
where the two left down arrows are norm maps and the right down
arrow is the multiplication by $[\F_{K'}:\F_K]$. The diagram is
commutative because for a divisor $D$ of $K'$ we have
$\deg_K(\mathbf{N}_{K'/K}(D))=[\F_{K'}:\F_K]\cdot\deg_{K'}(D)$.
Taking projective limit,  we get the exact sequence
\begin{equation}\label{e:w}
0\longrightarrow \Cl_{\K_{\infty},p}\longrightarrow
\W_{\K_{\infty},p} \stackrel{\dag}{\longrightarrow} \Z_p
\end{equation}
where $\dag$ is induced from the degree maps.

\begin{lemma}\label{l:dag}
The map $\dag$ is athe zero map if the constant $\Z_p$-extension
$\F_{p^{\infty}}k$ is contained in $\K_{\infty}$; otherwise the map
is surjective.
\end{lemma}
\begin{proof}
Suppose $K'$ and $K$ are the $m$'th and the $n$'th layers of
$\K_{\infty}$, with $m>n\gg0$, If the constant $\Z_p$-extension is
contained in $\K_{\infty}$ then we have $[\F_{K'}:\F_K]=p^{m-n}$.
Since $m$ can be arbitrary large, we must have $\dag=0$. On the
other hand, if the constant $\Z_p$-extension is not contained in
$\K_{\infty}$, then $[\F_{K'}:\F_K]=1$ and hence $\dag$ is
surjective.
\end{proof}

\begin{corollary}\label{c:charw}
The group $\W_{\K_{\infty},p}$ is a torsion
finitely generated module over $\Lambda_{\Upsilon}$. If $e=1 $ and $\K_{\infty}$ is not a
constant field extension, then
$$\chi_{{\Upsilon}}(\W_{\K_{\infty},p})=\chi_{{\Upsilon}}(\Cl_{\K_{\infty},p})
\cdot (s);$$ otherwise, we have
$$\chi_{{\Upsilon}}(\W_{\K_{\infty},p})=\chi_{{\Upsilon}}(\Cl_{\K_{\infty},p}).$$
\end{corollary}
\begin{proof}
Suppose we have chosen a basis of $\Upsilon$ and use it to identify
$\Lambda_{\Upsilon}$ with the formal power series ring
$\Z_p[[t_1,...,t_e]]$. Then the $\Lambda_{\Upsilon}$-module $\Z_p$
is isomorphic to $\Lambda_{\Upsilon}/(t_1,...,t_e)$, which is
pseudo-null unless $e=1$. If $e=1$, then $\Upsilon=\Xi$ and
$\chi_{{\Upsilon}}(\Z_p)=(s)$.
\end{proof}


To find the characteristic ideal of
$\W_{\K_{\infty},p}$ we consider a pseudo-isomorphism
\begin{equation}\label{e:psphi}
\phi:\W_{\K_{\infty},p}\longrightarrow \bigoplus_{j=1}^J
\Lambda_{\Upsilon}/(\zeta_j^{m_j})
\end{equation}
where each $\zeta_j$ is a prime element in $\Lambda_{\Upsilon}$.
Since any annihilator of a given element in
$\Lambda_{\Upsilon}/(\zeta_j^{m_j})$ must be inside the prime ideal
$(\zeta_j)$ which is of height one, the module
$\Lambda_{\Upsilon}/(\zeta_j^{m_j})$ contains no non-trivial
pseudo-null sub-module. Thus the kernel of $\phi$ is the maximal
pseudo-null submodule of $\W_{\K_{\infty},p}$. We use
$\N_{\K_{\infty}}$ to denote it.

We have chosen $s\in\Lambda_{\Xi}$ such that $\xi:=s+1$ is a topological
generator of $\Xi$.
Since $\N_{\K_{\infty}}$ is the maximal pseudo-null submodule, we
must have $\N_{\K_{\infty}}\cap
s\W_{\K_{\infty},p}=s\N_{\K_{\infty}}$, and hence
\begin{equation}\label{e:nsn}
\N_{\K_{\infty}}/s\N_{\K_{\infty}}\hookrightarrow
\W_{\K_{\infty},p}/s\W_{\K_{\infty},p}.
\end{equation}

\end{subsection}

\begin{subsection}{The $\Xi$-invariant part of $\W_{\K_{\infty},p}$}\label{subsec:inv}
For each open
subgroup $\alpha\subset \Upsilon$ let $\K_{\alpha}$ be the fixed
field of $\alpha$.
Then $\K_{p^n\Upsilon}$ is nothing but the $n$'th layer $\K_n$.
And the sub-system $\{\K_n\}_n$ is cofinal in the system
$\{\K_{\alpha}\}_{\alpha}$.  For an element
$x\in \W_{\K_{\infty},p}$, we use $x_{\alpha}\in\W_{\K_{\alpha},p}$
and $x_n\in\W_{K_n,p}$ to denote its images under the corresponding
natural maps.

\begin{lemma}\label{l:inv}
If $x\in \W_{\K_{\infty},p}^{\Xi}$, then for each open subgroup
$\alpha$ of $\Upsilon$ the divisor class $x_{\alpha}$ is
represented by a $\Xi$-invariant divisor of $\K_{\alpha}$.
\end{lemma}

\begin{proof}
For an open subgroup $\alpha$, let $D_{\alpha}$ be a $\Z_p$-divisor
of $\K_{\alpha}$ representing the class $x_{\alpha}$. Then there is
an $a_{\alpha}$ in the $p$-completion
${\widetilde {\K_{\alpha}^{\times}}}$   of
$\K_{\alpha}^{\times}$ such that
$(a_{\alpha})={}^{\xi}D_{\alpha}-D_{\alpha}$.
($\xi:=s+1$ is a topological
generator $\Xi$.)
This means that the image of $a_{\alpha}$ under the norm map from
$\K_{\alpha}^{\times}$ to $\K_{\alpha+\Xi}^{\times}$ gives rise to a trivial element of
${\widetilde { {\Prin}_{\K_{\alpha+\Xi}}    }}$. As we have observed that
${\widetilde {\K_{\alpha+\Xi}^{\times}}}$  and  ${\widetilde { {\Prin}_{\K_{\alpha+\Xi}}    }}$
are isomorphic
 this norm is the trivial element in ${\widetilde {\K_{\alpha+\Xi}^{\times}}}$ .
By Hilbert's Theorem 90 there is a element $b_{\alpha}$ in the
$p$-completion  ${\widetilde {\K_{\alpha}^{\times}}}$  such that
$(a_{\alpha})={}^{\xi}(b_{\alpha})-(b_{\alpha})$.
This shows that $x_{\alpha}$ is represented by the divisor
$D_{\alpha}-(b_{\alpha})$ which is invariant by the action of $\Xi$.
\end{proof}

\begin{lemma}\label{l:cs1}
If $x\in \W_{\K_{\infty},p}^{\Xi}$, then for each open subgroup
$\alpha$ of $ \Upsilon$ the divisor class $x_{\alpha}$ is
represented by a $\Xi$-invariant divisor which is supported on the set
$\cS_1(\K_{\alpha})$ consisting of places of $\K_{\alpha}$
sitting over $\cS_1$.
\end{lemma}

\begin{proof}
Choose $v_0\in\cS$ such that $v_0\in\cS_1$ if
$\cS_1\not=\emptyset$. Suppose that $v\in\cS$ is outside $\cS_1$
and $\cS':=\cS\setminus\{v\}$. Then we set a basis
$\sigma_1=\xi,\sigma_2,...,\sigma_e$ of $\Upsilon$ over $\Z_p$
such that $\Upsilon_v^0\subset \Z_p\sigma_2+\dots+\Z_p\sigma_e$.
For simplicity, denote
\begin{equation}\label{e:pair}
(m,n)=p^m\Z_p\sigma_1+p^n (\Z_p\sigma_2+\dots+\Z_p\sigma_e).
\end{equation}
The we have $\K_n=\K_{(n,n)}$. If $m \geq n$, then
$$\Gal(\K_{(m,n)}/\K_n)\simeq \frac{(n,n)}{(m,n)}\simeq
\frac{p^n\Z_p\cdot\xi}{p^m\Z_p\cdot\xi}$$
is a cyclic group generated by the restriction of $\xi$ on $\K_{(m,n)}$.
Since $\Upsilon_v^0\cap (n,n)=\Upsilon_v^0\cap (m,n)$,
the extension $\K_{(m,n)}/\K_n$ is unramified at every place sitting over $v$ and hence
unramified
outside $\cS'(\K_n)$.

There is a natural
embedding $i:\Div_{\K_n}\otimes_{\Z}\Z_p\hookrightarrow
\Div_{\K_{(m,n)}}\otimes_{\Z}\Z_p,$
sending a place (a prime divisor) $w_0$ of $\K_n$ to the divisor
$\sum_{w| w_0} e(w)w$ of $\K_{(m,n)}$, where $w$ runs through places
of $\K_{(m,n)}$ sitting over $w_0$ and $e(w)$ is the
ramification index of $w$ over $w_0$.
In particular, we have $e(w)=1$ if $w\not\in \cS_1(\K_{(m,n)})$.
Let $E=\sum_{w\not\in\cS_1(\K_{(m,n)})} a(w)w$ be a
$\Xi$-invariant $\Z_p$-divisor of $\K_{(m,n)}$  supported
outside $\cS_1(\K_{(m,n)})$. Then the action of
the Galois group $\Gal(\K_{(m,n)}/\K_n)$, which is the restriction
of the action of $\Xi$, fixes $E$ and we have $a(w)=a(w')$ for $w$,
$w'$ sitting over the same place of $\K_n$.
Therefore, $E$ is in
the image of the natural embedding $i$. Also, since $E$ is fixed by
the action of the Galois group $\Gal(\K_{(m,n)}/\K_n)$ the norm
$\mathbf{N}_{\K_{(m,n)}/\K_n}(E)$ just equals $p^{m-n}E$.

Let $D_{(m,n)}$ be a $\Xi$-invariant
$\Z_p$-divisor of $\K_{(m,n)}$ representing $x_{(m,n)}$, and put
$D_{(m,n)}=D^{(1)}_{(m,n)}+D^{(2)}_{(m,n)}$
where
$D^{(1)}_{(m,n)}$ is supported on $\cS'(\K_{(m,n)})$ and
$D^{(2)}_{(m,n)}$ is supported outside  $\cS'(\K_{(m,n)})$.
Both $D^{(1)}_{(m,n)}$ and $D^{(2)}_{(m,n)}$ are $\Xi$-invariant.
Put $E=D^{(2)}_{(m,n)}$, and from the above discussion we see that
if $D_n=\mathbf{N}_{\K_{(m,n)}/\K_n}(D_{(m,n)})$ then
\begin{equation}\label{e:ddd}
D_n=D^{(1)}+p^{m-n}E,
\end{equation}
where $D^{(1)}$ is a $\Xi$-invariant divisor of
$\K_n$ supported on $\cS'(\K_n)$.

Let $\iota,\varsigma\in\Z$ and $\kappa\in\Z_p^*$ be such that
$\deg (v_0)=\kappa\cdot p^{\iota}$ and $|\Cl_{\K_{(0,n)},p}|=p^{\varsigma}$.
These numbers are independent
of the choice of $m$. The  the divisor $\kappa p^{\iota}
E-\deg_{\K_{(0,n)}}(E)\cdot v_0$ is of degree zero, and hence its
multiple by $p^{\varsigma}$ is in the trivial divisor class. Thus,
there is an element $a$ in the $p$-completion ${\widetilde{\K_n^{\times}}}$ such
that
\begin{equation}\label{e:rhooo}
(a)=p^{\varsigma}(\kappa p^{\iota} E-\deg_{\K_{(0,n)}}(E)\cdot v_0).
\end{equation}
If $m$ is chosen to be greater than the integer $n+\iota+\varsigma$,
then from (\ref{e:ddd}) and (\ref{e:rhooo}) we find that the divisor
$D_n$ is $\Z_p$-linearly equivalent to a $\cS'(\K_n)$-supported divisor which is also
invariant under the action of $\Xi$.

We can replace $\cS$ by $\cS'$ and repeat the above argument, if necessary.
In the case where $\cS_1$ is non-empty, this will lead to the conclusion that $x_n$ is
represented by a $\Xi$-invariant divisor supported on $\cS_1(\K_n)$ and hence the proof is
completed. If $\cS_1=\emptyset$, then the above method shows that $x_n$ is
represented by a $\Xi$-invariant divisor supported on places sitting over $v_0$.
We then apply the above
argument again by taking $\cS=\{v_0\}$ and $v=v_0$.
This time in the equation (\ref{e:ddd}) the divisor $D^{(1)}$ is trivial and
$D_n$ is divisible by $p^{m-n}$. Since $m$ can be chosen arbitrary large,
the class $x_n$ is $p$-divisible in $\W_{\K_n}$, which is a finite $\Z_p$-module.
Therefore, we must have $x_n=0$.
The proof is completed.

\end{proof}

Using a similar method, we can make some further reduction.

\begin{lemma}\label{l:cs2}
If $x\in \W_{\K_{\infty},p}^{\Xi}$, then for each open subgroup
$\alpha\subset \Upsilon$ the divisor class $x_{\alpha}$ is
represented by a $\Xi$-invariant divisor which is supported on
$\cS_2(\K_{\alpha})$.
\end{lemma}
\begin{proof}
Choose $v_0$ in $\cS_1$ such that $v_0\in\cS_2$ if
$\cS_2\not=\emptyset$. For $v\in \cS_1\setminus\cS_2$, we choose a
basis $\sigma_1,...,\sigma_e$ of $\Upsilon$ such that
$\Upsilon_v^0\subset \Z_p\sigma_2+\dots+\Z_p\sigma_e$
and the index
$$[\Upsilon_v\cap (\Z_p\sigma_2+\dots+\Z_p\sigma_e):\Upsilon_v^0]< \infty.$$

Then under the natural projection
$$
\begin{array}{rcl}
\Upsilon & \longrightarrow &
\Upsilon/\Z_p\sigma_2+\dots+\Z_p\sigma_e\\
\sigma & \mapsto & {\bar \sigma}\\
\end{array}
$$
the image of $\Upsilon_v$ must be nontrivial. In particular, it is
generated by $p^{n'}\bar \sigma_1$ for some $n'$. Define the
subgroup $(m,n)$ of $\Upsilon$ as in the proof of the previous lemma.
Since we have $\Upsilon_v^0\cap (m,n)=\Upsilon_v^0\cap (0,n)$ the
extension $\K_{(m,n)}/\K_{(0,n)}$ is unramified at $v$.
Its
Galois group is generated by the restriction of $\bar \sigma_1$ to
$\K_{(m,n)}$ and the decomposition subgroup is generated by the
restriction of $p^{n'}\bar\sigma_1$. We also consider the field
extension $\K_{(m,n)}/\K_n$. This extension is unramified at $v$,
since $\K_{(0,n)}$ is a sub-field of
$\K_n$. If $m>n>n'$, then for
every place of $\K_n$ sitting over $v$ there is only one place of
$\K_{(m,n)}$ sitting over it. From this we make the key observation
that if $E$ is a divisor of $\K_{(m,n)}$ supported on places sitting
over $v$ then $E$ is in fact a divisor of $\K_n$. In this case, we
have
$\mathbf{N}_{\K_{(m,n)}/\K_n}(E)=p^{m-n}E$.

Now $x_{(m,n)}$ is represented by a $\Xi$-invariant divisor
$D_{(m,n)}$ in $\K_{(m,n)}$, which is supported on
$\cS_1(\K_{(m,n)})$. Let $E$ be the part of $D_{(m,n)}$ supported on
places sitting over $v$ and let $D^{(1)}_{(m,n)}=D_{(m,n)}-E$ which
is supported on $\cS_1'(\K_{(m,n)})$ where
$\cS_1':=\cS_1\setminus\{v\}$. Then both $E$ and $D^{(1)}_{(m,n)}$
are $\Xi$-invariant. Let $D^{(1)}$ be the image of $D^{(1)}_{(m,n)}$
under the norm map from $\K_{(m,n)}$ to $\K_n$. Then $D^{(1)}$ is a
divisor of $\K_n$ supported on $\cS_1'(\K_n)$ and
\begin{equation}\label{e:dddd}
D_n:=\mathbf{N}_{\K_{(m,n)}/\K_n}(D_{(m,n)})=D^{(1)}+p^{m-n}E.
\end{equation}
Then we finish the proof in the same way as the  last part of
the proof of Lemma \ref{l:cs1}.
\end{proof}

\end{subsection}

\begin{subsection}{Some special modules}\label{subs:special}
We will express the module $\W_{\K_{\infty}}^{\Xi}$ in terms of some
special modules. Suppose that $v\in\cS_1$. For each open subgroup
$\alpha\subset\Upsilon$, let $(\Upsilon/\alpha)_v$ be the
decomposition sub-group of $\Upsilon/\alpha=\Gal(\K_{\alpha}/k)$.
Put
\begin{equation}\label{e:wv}
\fG_{\Upsilon/\alpha,v}=\Lambda_{\Upsilon/\alpha}/\Im_{\Upsilon/\alpha,v},
\end{equation}
where $\Lambda_{\Upsilon/\alpha}=\Z_p[[\Upsilon/\alpha]]$ and
$\Im_{\Upsilon/\alpha,v}$ is the ideal generated by the set of all
$\sigma-1$ with $\sigma \in (\Upsilon/\alpha)_v+(\Xi+\alpha)/\alpha$.

For $\alpha\subset\beta$ we have the obvious homomorphism
$\rho_{\alpha,\beta,v}:\fG_{\Upsilon/\alpha,v}
\to \fG_{\Upsilon/\beta,v}$ which we use to form  the projective limit
\begin{equation}\label{e:wvim}
\fG_{v}=\lim_{\stackrel{\leftarrow}{\alpha}}\fG_{\Upsilon/\alpha,v}.
\end{equation}
We have
$
\fG_v=\Lambda_{\Upsilon}/\Im_v,
$
where the ideal $\Im_v$ is generated by the set
$\{\sigma-1\; | \; \sigma\in \Upsilon_v+\Xi\}$.

The group of $\Xi$-invariant divisors which are supported on
$\cS_2(\K_{\alpha})$ can be easily determined. For each $v\in\cS_2$,
we choose a place $v_{(\infty)}$ of $\K_{\infty}$ sitting over $v$,
and for each $\alpha$ let $v_{\alpha}$ be the place of $\K_{\alpha}$
sitting below $v_{(\infty)}$. The orbit $\T_{\alpha,v}$ of
$v_{\alpha}$ under the action of $\Xi$ is finite. Define the divisor
$D_{\alpha,v}=\sum_{w\in \T_{\alpha,v}} w $.
Under the action of $\Upsilon/\alpha$ the stabilizer of
$D_{\alpha,v}$ is the subgroup
$(\Upsilon/\alpha)_v+(\Xi+\alpha)/\alpha$. Then every
$\Xi$-invariant divisor of $\K_{\alpha}$ supported on places sitting
over $v$ can be express as ${\bar y}_vD_{\alpha,v}$ for some ${\bar
y}_v$ in the ring $\fG_{\Upsilon/\alpha}$. The assignment sending
$1$ to $[D_{\alpha,v}]$, the divisor class of $D_{\alpha,v}$, induces
a $\Lambda_{\Upsilon/\alpha}$-homomorphism
$\varphi_{\alpha,v}:\fG_{\Upsilon/\alpha,v}
\longrightarrow \W_{\K_{\alpha},p}^{\Xi}$. Taking  projective
limit, we get a $\Lambda_{\Upsilon}$-homomorphism
$\varphi_v:\fG_v\to  \W_{\K_{\infty},p}^{\Xi}$.

\begin{lemma}\label{l:invxi}
Then map
$$
\varphi=\sum_{v\in\cS_2}\varphi_v: \bigoplus_{v\in\cS_2}\fG_v
\longrightarrow  \W_{\K_{\infty},p}^{\Xi}
$$
is an isomorphism.
\end{lemma}

\begin{proof}
Lemma \ref{l:cs2} implies that $\varphi$ is surjective. Suppose that
${\bar y}=\sum_{v\in\cS_2}{\bar y}_v$ is in the kernel of $\varphi$.
Let $\rho_{\alpha,v}:\fG_v\longrightarrow \fG_{\Upsilon/\alpha,v}$
be the natural map, and denote the image $\rho_{\alpha,v}({\bar
y}_v)$ as ${\bar y}_{\alpha,v}$. We lift it through (\ref{e:wv}) to
an element $y_{\alpha,v}$ of the group ring
$\Lambda_{\Upsilon/\alpha}$. Then the divisor
$E_{\alpha}:=\sum_{v\in\cS_2} y_{\alpha,v}D_{\alpha,v}$ is
$\Z_p$-linear equivalent to zero.
We shall note that $E_{\alpha}$ is independent of the choice of
the above lifting.
To prove the injectivity of
$\varphi$, we only need to show that every $E_{\alpha}$ is the
trivial divisor, since this will imply that each $
y_{\alpha,v}D_{\alpha,v}$ is trivial and hence $ y_{\alpha,v}$ is in
the ideal $\Im_{\Upsilon/\alpha,v}$.

For $\alpha=p^n\Upsilon$ with $n$ large enough, each of the
intersections $\Upsilon_v\cap \alpha$, $v\in\cS_2$, contains
$\Xi\cap\alpha$ which is a direct
summand of $\alpha$. For a large $m$ let $\beta\subset \alpha$ be such that
the natural map $\Xi\cap\alpha\longrightarrow \alpha/\beta$ is surjective with
kernel equal $p^m(\Xi\cap\alpha)$. Then $\K_{\beta}/\K_{\alpha}$ is
a cyclic extension of degree $p^m$ and the decomposition subgroup at
each $v_{\alpha}$, $v\in\cS_2$, is the whole Galois group
$\Gal(\K_{\beta}/\K_{\alpha})$.
There are
$a_{\alpha}\in U_{\alpha}:=\O^{\times}_{\K_{\alpha},\cS_2}\otimes_{\Z}\Z_p$ and
$a_{\beta}\in U_{\beta}:=\O^{\times}_{\K_{\beta},\cS_2}\otimes_{\Z}\Z_p$ (the units
groups $\O^{\times}_{\K_{\alpha},\cS_2}$ and $\O^{\times}_{\K_{\beta},\cS_2}$ are
as those defined in Section \ref{subs:zeta}) such that
$E_{\alpha}=(a_{\alpha})$
and
$E_{\beta}=(a_{\beta})$.
Now $U_{\alpha}$ and
$U_{\beta}$ are of the same
rank over $\Z_p$. If $U_{\alpha}\varsubsetneq U_{\beta}$, then there
would be an $u\in \O^{\times}_{\K_{\beta},\cS_2} \setminus
\O^{\times}_{\K_{\alpha},\cS_2}$ such that $u_1:=u^p\in
\O^{\times}_{\K_{\alpha},\cS_2}$.
But this means that $u\in \K_{\alpha}(u_1^{\frac{1}{p}})\cap
\K_{\beta}=\K_{\alpha}$, since
$\K_{\alpha}(u_1^{\frac{1}{p}})/\K_{\alpha}$ is purely inseparable
while $\K_{\beta}/\K_{\alpha}$ is separable. This would lead to the
contradiction that $u\in \O^{\times}_{\K_{\alpha},\cS_2}$. Therefore
$U_{\alpha}=U_{\beta}$ and this implies that $E_{\beta}$ is a
divisor of $\K_{\alpha}$. Since $E_{\alpha}$ is the norm of the
divisor $E_{\beta}$ and hence equals $p^mE_{\beta}$, we conclude
that $E_{\alpha}$ and $a_{\alpha}$ are divisible by $p^m$. As $m$
can be arbitrary large, in the finite $\Z_p$-module $U_{\alpha}$ the
element $a_{\alpha}$ must be trivial. Therefore we have
$E_{\alpha}=0$.
\end{proof}
\begin{corollary}\label{c:invxi}
If no place in $\cS$ splits completely in  $\K_{\Xi}$, then
$\W_{\K_{\infty},p}^{\Xi}$ equals the kernel $\N_{\K_{\infty}}^{\Xi}$
of the pseudo-isomorphism $\phi$ {\em {(\ref{e:psphi})}}.
And we have an isomorphism $\N_{\K_{\infty}}^{\Xi}\simeq \oplus_{v\in\cS_2}\fG_v$.
\end{corollary}
\begin{proof}
Since for $v\in\cS$ the decomposition group $\Upsilon_v$ is not contained in
$\Xi$,
a prime ideal containing $\Im_v$ for some $v\in\cS_2$ must be of
height greater than one.
\end{proof}
\end{subsection}

\begin{subsection}
{The $\Xi$-co-invariant part of $\W_{\K_{\infty},p}$}
\label{subs:invariant}

Let $ \mathbf{N}_{\K_{\infty}/\K_{\Xi}}:\W_{\K_{\infty},p}\to
\W_{\K_{\Xi},p}$ be the norm map.
Let
$\mathcal{M_{\infty}}$ be the submodule $\ker(\mathbf{N}_{\K_{\infty}/\K_{\Xi}})/s\W_{\K_{\infty},p}$
of the quotient $\W_{\K_{\infty},p}/s\W_{\K_{\infty},p}$.
Define $\fG_{v}$ as in (\ref{e:wvim}).
Our goal is to establish the  exact sequence of Lemma \ref{l:summerize}
$$0\longrightarrow \mathcal{M}_{\infty}\longrightarrow
\bigoplus_{v\in\cS_1}\fG_v\longrightarrow \Z_p\longrightarrow 0.$$

For a topological group $G$ let  $G^{\wedge}$ denote the Pontryagin
dual group $\Hom_{cont}(G,\Q_p/\Z_p)$
consisting of continuous
homomorphisms.
As the duality in the case of the dual ${ \W}_{\K_{\infty},p}^{\wedge}$
of $\W_{\K_{\infty},p}$ respects the actions of $\Upsilon$,
it is a duality between $\Lambda_{\Upsilon}$-modules.
In particular, an element $\phi\in {\W}_{\K_{\infty},p}^{\wedge}$
annihilates $s\W_{\K_{\infty},p}$  if and only
if for all $x\in
\W_{\K_{\infty},p}$ we have
$$0=\phi((\xi-1)x)=({\xi^{-1}}\phi)(x)-\phi(x)=({\xi^{-1}}\phi-\phi)(x)$$
($\xi=s+1$ is the chosen topological generator of $\Xi$).
Therefore the annihilators of
$s\W_{\K_{\infty},p}$ are the elements of
$({\W}_{\K_{\infty},p}^{\wedge})^{{}^{\Xi}}$.
And the Pontryagin dual of the
quotient $\W_{\K_{\infty},p}/s\W_{\K_{\infty},p}$ is the
$\Xi$-invariant sub-module
$({\W}_{\K_{\infty},p}^{\wedge})^{{}^{\Xi}}$. Using this, we deduce in a similar way that
\begin{equation}\label{e:md}
\mathcal{M}_{\infty}^{\wedge}=
({\W}_{\K_{\infty},p}^{\wedge})^{{}^{\Xi}}/i^{*}(\W_{\K_{\Xi},p}^{\wedge}).
\end{equation}

Since $\W_{\K_{\infty},p}$ is identified with the Galois group of
the maximal everywhere unramified pro-$p$ abelian extension  over
$\K_{\infty}$, the module ${ \W}_{\K_{\infty},p}^{\wedge}$ is
identified as a sub-module of
$\Gal({\bar k}/\K_{\infty})^{\wedge}$ where ${\bar k}$
is a fixed separable closure of $k$.

Since ${ \W}_{\K_{\infty},p}$ is compact, the image of an
$\omega\in { \W}_{\K_{\infty},p}^{\wedge}$ is a finite, and hence cyclic,
subgroup of $\Q_p/\Z_p$. Therefore, the fixed field of $\ker(\omega)$
is a finite cyclic extension $\K_{\infty}^{(\omega)}$ over
$\K_{\infty}$.
If
$\omega$ in $({ \W}_{\K_{\infty},p}^{\wedge})^{{}^{\Xi}}$,
then the extension $\K_{\infty}^{(\omega)}/\K_{\infty}$
is invariant under the action of the
$\xi$. This implies that the extension $\K_{\infty}^{(\omega)}/\K_{\Xi}$ is also
abelian and the Galois group $G:=\Gal(\K_{\infty}^{(\omega)}/\K_{\Xi})$
is an extension of the Galois group
$\Gal(\K_{\infty}/\K_{\Xi})=\Xi\simeq \Z_p$ by the finite cyclic group
$\Gal(\K_{\infty}^{(\omega)}/\K_{\infty})$. Therefore, $G$ is the
direct product of a subgroup $H\simeq\Z_p$ with the finite $p$-torsion subgroup
$\Gal(\K_{\infty}^{(\omega)}/\K_{\infty})$.
Let $\K'$ be the fixed field of $H$. Then
$\K_{\infty}^{(\omega)}=\K_{\infty}\K'$ and
$\Gal(\K_{\infty}^{(\omega)}/\K_{\infty})\simeq \Gal(\K'/\K_{\Xi})$.
Thus, if we identify these two Galois groups, then $\omega$ can be
obtained form a character of $\Gal(\K'/\K_{\Xi})$.

Let
$i^{*} : \Gal({\bar k}/\K_{\Xi})^{\wedge} \to \Gal({\bar k}/\K_{\infty})^{\wedge}$
be the homomorphism dual to the
restriction map of the Galois groups.
We have just shown that
$({\W}_{\K_{\infty},p}^{\wedge})^{{}^{\Xi}}$
is contained in the image of $i^{*}$. Let
$F_{\infty}:={i^{*}}^{-1}(({\W}_{\K_{\infty},p}^{\wedge})^{{}^{\Xi}})$.
Let $\L_{\Xi}/\K_{\Xi}$ be the maximal everywhere unramified pro-$p$
abelian extension and let
$\fG'_{\K_{\Xi},p}=\Gal(\K_{\infty}\L_{\Xi}/\K_{\Xi})$.
Put $H_{\infty}:=(\fG'_{\K_{\Xi},p})^{\wedge}$.
Then
$F_{\infty}$ contains $H_{\infty}$ and the map $i^{*}$ induces an isomorphism
$F_{\infty}/H_{\infty}\simeq ({
\W}_{\K_{\infty},p}^{\wedge})^{{}^{\Xi}}/i^{*}(\W_{\K_{\Xi},p}^{\wedge})$.
This and (\ref{e:md})
imply the following lemma.

\begin{lemma}\label{l:mfh}
The $\Lambda_{\Upsilon}$-module $\mathcal{M}_{\infty}$ is dual to the quotient
$F_{\infty}/H_{\infty}$.
\end{lemma}

Let $\alpha$  be an open
subgroup $\alpha$  of $\Upsilon$ containing $\Xi$. Write $\K_{\alpha}$  for the fixed field of ${\alpha}.$
Let $\A_{\K_{\alpha}}^{\times}$ be the idele group of $\K_{\alpha}$ and let
$$\psi_{\alpha}:\A_{\K_{\alpha}}^{\times}\longrightarrow \Gal(\K_{\infty}/\K_{\alpha})=\alpha$$
be the global norm residue map. Then for each place $w$ of
$\K_{\alpha}$ the local norm residue map can be viewed as the composition
$$\psi_{w,\alpha}:\K_{\alpha,w}^{\times}\longrightarrow
\A_{\K_{\alpha}}^{\times}\stackrel{\psi_{\alpha}}{\longrightarrow}
\alpha.$$ Let $\mathcal{C}_w\subset \O_{w}^{\times}$ be the intersection
$\ker(\psi_{w,\alpha})\cap\O_w^{\times}$. Then through the local norm residue map
$\O_w^{\times}/\mathcal{C}_w$ is identified with the inertia subgroup
$\alpha_w^0$ of $\alpha$. Denote $\daleth_w=\Xi\cap \alpha_w^0$ and
put $\mathcal{B}_w=\psi_{w,\alpha}^{-1}(\daleth_w)$. Then we have
the exact sequence
\begin{equation}\label{e:bc}
0\longrightarrow \mathcal{B}_w/\mathcal{C}_w\stackrel{{\bar
\psi}_{w,\alpha}} {\twoheadrightarrow} \daleth_w\subset \Xi.
\end{equation}

Let
$\mathcal{G}_{\alpha}$  be the $p$-completion
of the idele class group
$$\K_{\alpha}^{\times}\backslash \A_{\K_{\alpha}}^{\times}/
\prod_{w\in\cS_1(\K_{\alpha})}\mathcal{C}_w\cdot
\prod_{w\not\in\cS_1(\K_{\alpha})} \O_w^{\times}$$
and let  ${\bar {\mathcal{G}}}_{\alpha}$ be that
of the idele class group
$$\K_{\alpha}^{\times}\backslash \A_{\K_{\alpha}}^{\times}/
\prod_{w\in\cS_1(\K_{\alpha})}\mathcal{B}_w\cdot
\prod_{w\not\in\cS_1(\K_{\alpha})} \O_w^{\times}.$$
It is not difficult to see that the kernel of the natural map
$Q:\mathcal{G}_{\alpha}\longrightarrow {\bar
{\mathcal{G}}}_{\alpha}$
is exactly the image of the natural embedding
$\prod_{w\in\cS_1(\K_{\alpha})}\mathcal{B}_w/\mathcal{C}_w
\longrightarrow
\mathcal{G}_{\alpha}$.
And in view of (\ref{e:bc})
it is obvious that $\psi_{\alpha}(\ker(Q))\subset\Xi$.
In other words, we can define a map
$\gimel_{\alpha}:\prod_{w\in\cS_1(\K_{\alpha})}\mathcal{B}_w/\mathcal{C}_w
\longrightarrow \Xi$ and incorporate these in the following
commutative diagram:
$$
\begin{array}{ccc}
\mathcal{G}_{\alpha} & \stackrel{\psi_{\alpha}}{\longrightarrow} & \alpha\\
\uparrow & \circlearrowleft & \cup \\
\prod_{w\in\cS_1(\K_{\alpha})}\mathcal{B}_w/\mathcal{C}_w &
\stackrel{\gimel_{\alpha}}{\longrightarrow} & \Xi.\\
\end{array}
$$
The next lemma involves the projective limit of modules of the form
$\mathcal{A}_{\alpha}:=\ker(\gimel_{\alpha})$.
If
$\Xi\subset\beta\subset \alpha$ then the norm map on ideles induces
a surjective homomorphism
$\mathcal{A}_{\beta}\twoheadrightarrow\mathcal{A}_{\alpha}$ and we
 denote the projective limit by:
$\mathcal{A}_{\infty}=\lim_{\stackrel{\leftarrow}{\alpha}}\mathcal{A}_{\alpha}$.
It is easy to see that
$\mathcal{A}_{\alpha}=\ker(Q)\cap\ker(\psi_{\alpha})$.

\begin{lemma}\label{l:aalpha}
The $\Lambda_{\Upsilon}$-module $F_{\infty}/H_{\infty}$  is dual to $\mathcal{A}_{\infty}$ .
\end{lemma}
\begin{proof}
Firstly we express  $F_{\infty}$ as a direct limit. For
each open sub-group $\alpha$ of $\Upsilon$ containing $\Xi$,  we define
$F_{\alpha}$  as the subgroup of
the group
$\Gal({\bar k}/\K_{\alpha})^{\wedge}$ consisting of
elements $\omega$ such that $j_{\alpha}^{*}(\omega)\in F_{\infty}$
, where
$ j_{\alpha}^{*} :  \Gal({\bar k}/\K_{\alpha})^{\wedge} \to \Gal({\bar k}/\K_{\Xi})^{\wedge}$
is the dual of the
restriction of Galois group. Since every cyclic extension over
$\K_{\Xi}$ is obtained from some cyclic extension over $\K_{\alpha}$
for some $\alpha$, the modules $F_{\infty}$  is a
direct limit: $F_{\infty}=\lim_{\stackrel{\rightarrow}{\alpha}}F_{\alpha}$
where $\alpha$ runs through all the open
subgroups of $\Upsilon$ containing $\Xi$.
Let $\K'/\K_{\alpha}$ be an abelian extension with
$\Gal(\K'/\K_{\alpha})=\beth$, and let
$\psi:\A_{\K_{\alpha}}^{\times}\longrightarrow\beth$ and
$\psi_w:\K_w^{\times}\longrightarrow \beth_w$ be the corresponding global
and local norm residue maps. Under this setting, the condition
$\K_{\infty}\subset\K'$ is equivalent to $\ker(\psi)\subset
\ker(\psi_{\alpha})$. If this holds, then the condition that
$\K'/\K_{\infty}$ is unramified at places sitting over $w$ is
equivalent to the condition that $\mathcal{C}_w\subset\ker(\psi_w)$.
Therefore $\mathcal{G}_{\alpha}$ is the Galois group over
$\K_{\alpha}$ of the maximal pro-$p$ abelian extension containing
$\K_{\infty}$ such that it is everywhere unramified over
$\K_{\infty}$. It is obvious that ${
{\mathcal{G}}}_{\alpha}^{\wedge}\subset F_{\alpha}$. If $\omega\in
F_{\alpha}$ and $\K^{(\omega)}_{\alpha}$ is the fixed field of its
kernel, then the extension
$\K_{\infty}\K^{(\omega)}_{\alpha}/\K_{\alpha}$ is abelian and the
extension $\K_{\infty}\K^{(\omega)}_{\alpha}/\K_{\infty}$ is
everywhere unramified. This means that
$\omega\in { {\mathcal{G}}}_{\alpha}^{\wedge}$. Hence we see that
$F_{\alpha}={ {\mathcal{G}}}_{\alpha}^{\wedge}$.

Secondly we consider $H_{\infty}$.
As before  for an open subgroup $\alpha$ of $\Upsilon$ containing $\Xi$
we put $ H_{\alpha}$  as the subgroup of
the group $\Gal({\bar k}/\K_{\alpha})^{\wedge}$
consisting of
elements $\omega$ such that $j_{\alpha}^{*}(\omega)\in H_{\infty}$. Then we have
$H_{\infty}=\lim_{\stackrel{\rightarrow}{\alpha}}H_{\alpha}$.
Again, let $\K'/\K_{\alpha}$ be an abelian extension with
$\Gal(\K'/\K_{\alpha})=\beth$ and let
$\psi:\A_{\K_{\alpha}}^{\times}\longrightarrow\beth$ and
$\psi_w:\K_{\alpha,w}^{\times}\longrightarrow \beth_w$ be the corresponding
global and local norm residue maps. The condition that $\K'$
contains the field $\K_{\Xi}$ is equivalent to the condition that
$\ker(\psi)\subset \psi_{\alpha}^{-1}(\Xi)$. If this holds, then the
condition that $\K'/\K_{\Xi}$ is unramified at places sitting over
$w$ is equivalent to the condition that
$\mathcal{B}_w\subset\ker(\psi_w)$. This implies that ${\bar
{\mathcal{G}}}_{\alpha}$ is the Galois group over $\K_{\alpha}$ of
the maximal pro-$p$ abelian extension containing $\K_{\Xi}$ such
that it is everywhere unramified over $\K_{\Xi}$. Denote this field
extension as $\L'/\K_{\alpha}$. Then $H_{\alpha}$ is dual to the Galois
group $\Gal(\K_{\infty}\L'/\K_{\alpha})$. But since $k_{\infty}$ is the fixed field
of $\ker (\psi_{\alpha})\subset \mathcal{G}_{\alpha}$ and $\L'$ is that of
$\ker(Q)\subset \mathcal{G}_{\alpha}$
we have the natural
isomorphism
$\Gal(\K_{\infty}\L'/\K_{\alpha})=\mathcal{G}_{\alpha}/\ker(Q)\cap\ker(\psi_{\alpha})$.
This shows that $H_{\alpha}=(\mathcal{G}_{\alpha}/\ker(Q)\cap\ker(\psi_{\alpha}))^{\wedge} $.
In particular, we see that inside $F_{\alpha}=G_{\alpha}^{\wedge}$
the subgroup $H_{\alpha}$ is the annihilator of the subgroup
$\ker(Q)\cap\ker(\psi_{\alpha})\subset G_{\alpha}$, and hence
$F_{\alpha}/H_{\alpha}=(\ker(Q)\cap\ker(\psi_{\alpha}))^{\wedge}$.
But we have seen that
$\mathcal{A}_{\alpha}$ equals to the
intersection $\ker(Q)\cap\ker(\psi_{\alpha})$.
The proof of the lemma is completed.
\end{proof}

Finally we relate $\mathcal{A}_{\infty}$ to the special modules $\fG_v$ of section \ref{subs:special}.

\begin{lemma}
We have  an exact sequence of $\Lambda_{\Upsilon}$-modules
$$
0 \longrightarrow \mathcal{A}_{\infty}\longrightarrow
\bigoplus_{v\in\cS_1}\fG_v\longrightarrow \Z_p\longrightarrow 0,
$$
where  $\Z_p$ is endowed
with the trivial action of $\Upsilon$.
\end{lemma}
\begin{proof}
We
first note that if $w\not\in \cS_1(\K_{\alpha})$ then the group
$\daleth_w$ is trivial and $\mathcal{B}_w=\mathcal{C}_w$; otherwise
$\daleth_w$ is non-trivial and hence isomorphic to $\Z_p$. And we
also observe that the group $\daleth_w$ only depends on the place $v\in \cS_1$ sitting below
$w$. In fact, if $\Xi\subset\beta\subset\alpha$ and $w'$ is a
place of $\K_{\beta}$ sitting over $v$, from the definition, we see
that
\begin{equation}\label{e:ww'}
\daleth_{w'}= \daleth_w.
\end{equation}
For simplicity, we denote $\daleth_v=\daleth_w$. To treat the groups
$\daleth_v$, $v\in\cS_1$ in a consistent way, we let $d_v$ denote the
integer such that $p^{d_v}\xi$ is a generator of $\daleth_v$.
As before, for each place $v\in\cS_1$, choose a place $v_{(\infty)}$
of $\K_{\infty}$ sitting over $v$ and for each $\alpha$
denote by $v_{\alpha}$ the
place of $\K_{\alpha}$ sitting below $v_{(\infty)}$. Every place
$w\in\cS_1(\K_{\alpha})$ sitting over $v$ is in the orbit of
$v_{\alpha}$ under the action of $\Upsilon$. Thus, there is a
$\sigma\in\Upsilon$ such that $w=\sigma(v_{\alpha})$. In this case,
we have the commutative diagram:
\begin{equation}\label{e:com}
\begin{array}{ccc}
\mathcal{B}_{v_{\alpha}}/\mathcal{C}_{v_{\alpha}} & \stackrel{{\bar
\psi}_{v_{\alpha},\alpha}}
{\longrightarrow} & \daleth_{v}\\
\downarrow \sigma & \circlearrowright & ||\\
\mathcal{B}_w/\mathcal{C}_w &  \stackrel{{\bar \psi}_{w,\alpha}}{\longrightarrow} & \daleth_{v}.\\
\end{array}
\end{equation}
Recall that the homomorphism ${\bar \psi}_{v_{\alpha,\alpha}}$ is the one
in (\ref{e:bc}). Put
$b_{v,\alpha}={\bar \psi}_{v_{\alpha,\alpha}}^{-1}(p^{d_v}\xi)$.
Since under the action of $\Upsilon$, the stabilizer of $v_{\alpha}$
is $\Upsilon_{v}+\alpha$, where it is assumed that
$\Xi\subset\alpha$, we have an  isomorphism of
$\Upsilon$-modules:
$\digamma_{v,\alpha}:\fG_{\Upsilon/\alpha,v}
{\simeq}  \prod_{w\mid v} \mathcal{B}_w/\mathcal{C}_w
$
given by $z_v   \mapsto    z_v b_{v,\alpha}$,
and a  commutative diagram
\begin{equation}\label{e:digammamcom}
\begin{array}{ccl}
\bigoplus_{v\in\cS_1}\fG_{\Upsilon/\alpha,v} & \stackrel{\Sigma_{\alpha}}{\longrightarrow} & \Z_p\\
\downarrow  & \circlearrowright & \downarrow\\
\prod_{w\in\cS_1(\K_{\alpha})}\mathcal{B}_w/\mathcal{C}_w &
\stackrel{\gimel_{\alpha}}{\longrightarrow} & \Xi\\
\end{array}
\end{equation}
where the ring $\fG_{\Upsilon/\alpha,v}$ is defined in Section
\ref{subs:special}, the left down-arrow is the map
$\digamma_{\alpha}$ sending $\sum_{v}z_v$ to $\sum_v
\digamma_{v,\alpha}(z_v)$, the right down-arrow sends $1$ to $\xi$
and the map $\Sigma_{\alpha}$ sends each unit $1_{\alpha,v}$ in the
ring $\fG_{\Upsilon/\alpha,v}$ to $p^{d_v}\xi$. Therefore, we have
an isomorphism
\begin{equation}\label{e:kerdi}
\digamma_{\alpha}{}_{\mid_{\ker(\Sigma_{\alpha})}}:\ker(\Sigma_{\alpha})
\stackrel{\sim}{\longrightarrow} \mathcal{A}_{\alpha}.
\end{equation}
Now we consider the projective limits of the above objects. First it
is easy to see that
$$\bigoplus_{v\in\cS_1}\fG_v
=\lim_{\stackrel{\leftarrow}{\alpha}}\bigoplus_{v\in\cS_1}\fG_{\Upsilon/\alpha,v}$$
where
$\alpha$ runs through all open sub-group of $\Upsilon$ containing
$\Xi$. Since $\bigoplus_{v\in\cS_1}\fG_v$ is compact and the image
of $\Sigma_{\alpha}$ is an open sub-group of $\Z_p$ independent of
$\alpha$, we have an exact sequence
$0\longrightarrow \lim_{\stackrel{\leftarrow}{\alpha}}\ker(\Sigma_{\alpha})
\longrightarrow \bigoplus_{v\in\cS_1}\fG_v \longrightarrow
\Z_p\longrightarrow 0$.
 If $\Xi\subset\beta\subset\alpha$ and $w'$
is a place of $\K_{\beta}$ sitting over $w$, then we have the
commutative diagram
\begin{equation}\label{e:comcom}
\begin{array}{rcl}
\mathcal{B}_{w'}/\mathcal{C}_{w'} & \stackrel{{\bar \psi}_{w',\beta}}{\longrightarrow} & \daleth_{v}\\
\downarrow  & \circlearrowright & ||\\
\mathcal{B}_w/\mathcal{C}_w &  \stackrel{{\bar \psi}_{w,\alpha}}{\longrightarrow} & \daleth_{v},\\
\end{array}
\end{equation}
where the left down-arrow is from the local norm map. This together
with the isomorphism (\ref{e:kerdi}) and the diagram
(\ref{e:digammamcom}) implies
$\mathcal{A}_{\infty}\simeq
 \lim_{\stackrel{\leftarrow}{\alpha}}\ker(\Sigma_{\alpha})$.
\end{proof}

We summarize the above discussions in the following lemma.

\begin{lemma}\label{l:summerize}
Let
$$\mathbf{N}:\W_{\K_{\infty},p}/s\W_{\K_{\infty},p}\longrightarrow
\W_{\K_{\Xi},p}$$ be the natural map induced from the norm
$\mathbf{N}_{\K_{\infty}/\K_{\Xi}}$. If $\K_{\infty}/\K_{\Xi}$ is
everywhere unramified, then we have the exact sequence of
$\Lambda_{\Upsilon}$-modules
$$ 0\longrightarrow \W_{\K_{\infty},p}/s\W_{\K_{\infty},p}\stackrel{\mathbf{N}}{\longrightarrow}
\W_{\K_{\Xi},p} \longrightarrow \Z_p\longrightarrow 0;$$ otherwise,
we have exact sequences of $\Lambda_{\Upsilon}$-modules
$$0\longrightarrow \mathcal{M}_{\infty}\longrightarrow
\W_{\K_{\infty},p}/s\W_{\K_{\infty},p}\stackrel{\mathbf{N}}{\longrightarrow}
\W_{\K_{\Xi},p}\longrightarrow \mathcal{Z}\longrightarrow 0
$$
and
$$0\longrightarrow \mathcal{M}_{\infty}\longrightarrow
\bigoplus_{v\in\cS_1}\fG_v\longrightarrow \Z_p\longrightarrow 0,$$
where $\mathcal{Z}$ is of finite cardinality and $\Z_p$ is endowed
with the trivial action of $\Upsilon$.
\end{lemma}

\begin{proof}
If $\K_{\infty}/\K_{\Xi}$ is everywhere unramified, then $\cS_1$ is
empty and $\mathcal{M}_{\infty}$ is trivial. We need to determine
the cokernel of $\mathbf{N}$. From the duality, we see that it is
dual to the kernel of
$i^{*}\mid_{{ \W}_{\K_{\Xi},p}^{\wedge}}:{\W}_{\K_{\Xi},p}^{\wedge} \to
\Gal({\bar k}/\K_{\infty})^{\wedge}.
$
It is easy to see that this kernel is just the intersection
${\Xi}^{\wedge} \cap { \W}_{\K_{\Xi},p}^{\wedge}$. Here we consider
both ${\Xi}^{\wedge}$ and ${ \W}_{\K_{\Xi},p}^{\wedge}$ as subgroups
of $\Gal({\bar k}/\K_{\Xi})^{\wedge}$. Since
$\W_{\K_{\Xi},p}$ is the Galois group of the maximal everywhere
unramified pro-$p$ abelian extension of $\K_{\Xi}$, the intersection
is the Pontryagin dual of the quotient $\Xi/\Xi'$ where $\Xi'\subset
\Xi$ is the subgroup generated by all the inertia groups at all the
places. Consequently, if $\K_{\infty}/\K_{\alpha}$ is everywhere
unramified then $\coker(\mathbf{N})\simeq\Z_p$; otherwise, it is a
finite set.
\end{proof}

\begin{corollary}\label{c:summerize}
If $e>1$ and no place in $\cS_1$ splits completely under
$\K_{\Xi}/k$, then the $\Lambda_{\Upsilon}$-modules
$\mathcal{M}_{\infty}$ and $\W_{\K_{\infty},p}/s\W_{\K_{\infty},p}$
are pseudo-null $\Lambda_{\Upsilon}$-modules.
\end{corollary}

\begin{proof}
Since $\Xi\varsubsetneq\Upsilon$, a finitely generated torsion
module over $\Lambda_{\Upsilon/\Xi}$ is pseudo-null if it is
considered as a $\Lambda_{\Upsilon}$-module. Since $\Z_p$ and
$\W_{\K_{\Xi},p}$ are finitely generated torsion over
$\Lambda_{\Upsilon/\Xi}$, in view of Lemma \ref{l:summerize}, we
only need to show that for each $v\in\cS_1$, the module $\fG_v$ is
also finitely generated torsion over $\Lambda_{\Upsilon/\Xi}$. But
we have
$$\fG_v=\Lambda_{\Upsilon}/\Im_v\simeq \Lambda_{\Upsilon/\Xi}/{\bar \Im}_v$$
where ${\bar \Im}_v$ is the ideal of $\Lambda_{\Upsilon/\Xi}$
generated by the set $\{\sigma-1\; | \; \sigma\in
\Upsilon_v+\Xi/\Xi\}$. Since $v$ does not completely split over
$\K_{\Xi}/k$, the quotient group $\Upsilon_v+\Xi/\Xi$ is non-trivial
and hence $\fG_v$ is torsion over $ \Lambda_{\Upsilon/\Xi}$.
\end{proof}
\end{subsection}

\begin{subsection}{Greenberg's lemma}\label{l:g}
\begin{lemma}{\em (Greenberg)}\label{l:grn}
Let $\Upsilon\simeq\Z_p^e$ for some $e$ and let $\mathcal{Y}$ be a
finitely generated torsion $\Lambda_{\Upsilon}$-module. Then the
following are true.
\begin{enumerate}
\item Assume that $\mathcal{Y}$ has an annihilator $\Phi\in\Lambda_{\Upsilon}$ such that
$p\nmid \Phi$. Then $\Upsilon$ contains at least one sub-group
$\Upsilon'$ such that $\Upsilon/\Upsilon'\simeq\Z_p$ with the
property that $\Y$ is finitely generated over $\Lambda_{\Upsilon'}$.
Furthermore, for every $e-1$-dimensional $\F_p$-subspace ${\bar
\Upsilon}'$ of ${\bar \Upsilon}:=\Upsilon\otimes \Z_p/p\Z_p$, the
sub-group $\Upsilon'$ can be chosen such that its image under
$\Upsilon\longrightarrow {\bar \Upsilon}$ equals ${\bar \Upsilon}'$.
\item If $\Y$ is pseudo-null, then {\em (1)} holds. In this case $\Y$ is a torsion module over
$\Lambda_{\Upsilon'}$.
\end{enumerate}
\end{lemma}

\begin{proof}
The first part of statement (1) is actually Lemma 2 in \cite{grn},
and its proof actually proves the last part of (1). Statement (2) is
proved in the discussion after the proof of Lemma 2.
\end{proof}
\end{subsection}

\begin{subsection}{The characteristic ideals}\label{subs:ch}
In this subsection we shall compare the characteristic ideals
$\chi_{{\Upsilon}}(\Cl_{\K_{\infty},p})$ and
$\chi_{{\Upsilon/\Xi}}(\Cl_{\K_{\Xi},p})$ .
We assume that no place in $\cS$ splits completely
under $\K_{\Xi}/k$. As before we  let
$\phi:\W_{\K_{\infty},p} \to \oplus_{j=1}^J  \Lambda_{\Upsilon} /(\zeta_j^{m_j})$
be a pseudo-isomorphism with $\ker(\phi)=\mathcal{N}_{\K_{\infty}}$.
We also denote
$\mathcal{T}_{\K_{\infty}}=\coker(\phi)$.
Then we have exact sequences
\begin{equation}\label{e:nwim}
0\longrightarrow \mathcal{N}_{\K_{\infty}}\longrightarrow
\W_{\K_{\infty},p} \longrightarrow \image(\phi)\longrightarrow 0,
\end{equation}
and
\begin{equation}\label{e:imt}
0\longrightarrow \image(\phi)\longrightarrow \bigoplus_{j=1}^J
\Lambda_{\Upsilon}/(\zeta_j^{m_j}) \longrightarrow
\mathcal{T}_{\K_{\infty}}\longrightarrow 0.
\end{equation}
By Corollary \ref{c:invxi} and
Corollary \ref{c:summerize} the $\Lambda_{\Upsilon}$-module
$\Y=\N_{\K_{\infty}}\oplus \mathcal{T}_{\K_{\infty}}\oplus
\W_{\K_{\infty},p}/s\W_{\K_{\infty},p}\oplus \W_{\K_{\infty},p}^{\Xi}$ is pseudo-null. Applying
Lemma \ref{l:grn}, we can find a subgroup $\Upsilon'$ such that $\Y$
is a finitely generated torsion $\Lambda_{\Upsilon'}$-module. We can
choose $\Upsilon'$ such that the subspaces ${\bar \Upsilon}'$ and
${\bar \Xi}:=\Xi\otimes\Z_p/p\Z_p$ span the space $\bar\Upsilon$.
This means that $\Upsilon$ is the direct sum of $\Xi$ and
$\Upsilon'$.

The action of $s=\xi-1$ commutes with (\ref{e:nwim}) and
(\ref{e:imt}) and from the snake lemma, we have exact sequences
\begin{equation}\label{e:snake1}
\begin{array}{rc}
0\longrightarrow  \mathcal{N}_{\K_{\infty}}^{\Xi}\longrightarrow
\W_{\K_{\infty},p}^{\Xi} \longrightarrow
\image(\phi)^{\Xi}\longrightarrow &
\mathcal{N}_{\K_{\infty}}/s \mathcal{N}_{\K_{\infty}}\\
{} & \downarrow\\
0\longleftarrow \image(\phi)/s\image(\phi) \longleftarrow &
\W_{\K_{\infty},p}
 /s\W_{\K_{\infty},p}
\\
\end{array}
\end{equation}
and
\begin{equation}\label{e:snake2}
\begin{array}{rc}
0\longrightarrow  \image(\phi)^{\Xi}\longrightarrow
(\bigoplus_{i=1}^J \Lambda_{\Upsilon}/(\zeta_j^{m_j}))^{\Xi}
\longrightarrow  & \mathcal{T}_{\K_{\infty}}^{\Xi} \\
{} & \downarrow\\
0\longleftarrow  \mathcal{T}_{\K_{\infty}}/s
\mathcal{T}_{\K_{\infty}} \longleftarrow  \bigoplus_{j=1}^J
\Lambda_{\Upsilon/\Xi}/({\bar \zeta}_j^{m_j}) \longleftarrow &
\image(\phi)/s\image(\phi),
\\
\end{array}
\end{equation}
where ${\bar \zeta}_i$ is the image of $\zeta_i$ under the
projection
$\Lambda_{\Upsilon}\longrightarrow \Lambda_{\Upsilon/\Xi}.$

\begin{lemma}\label{l:nonzero}
Under the condition that no place in $\cS$ splits completely over
$\K_{\Xi}/k$, we have ${\bar \zeta}_i\not=0$, for each $j  \in J$.
\end{lemma}
\begin{proof}
If ${\bar \zeta}_i=0$ for some $i$, then
$\Lambda_{\Upsilon/\Xi}/({\bar \zeta}_i^{m_i})$ is a free
$\Lambda_{\Upsilon/\Xi}$-module. Since $\mathcal{T}_{\K_{\infty}}/s
\mathcal{T}_{\K_{\infty}}$, $\mathcal{T}_{\K_{\infty}}^{\Xi}$ and
$\mathcal{N}_{\K_{\infty}}/s \mathcal{N}_{\K_{\infty}}$ are all
torsion $\Lambda_{\Upsilon/\Xi}$-modules, the above exact sequences
say that neither $\image(\phi)/s\image(\phi)$ nor
$\W_{\K_{\infty},p}/s\W_{\K_{\infty},p}$ is a torsion
$\Lambda_{\Upsilon/\Xi}$-module. The homomorphism
$\Upsilon'\longrightarrow \Upsilon\longrightarrow \Upsilon/\Xi$
induces an identification of $\Lambda_{\Upsilon'}$ with
$\Lambda_{\Upsilon/\Xi}$. This implies that
$\W_{\K_{\infty},p}/s\W_{\K_{\infty},p}$ is not a torsion
$\Lambda_{\Upsilon'}$-module. On the other hand, we know that $\Y$
is a torsion $\Lambda_{\Upsilon'}$-module and hence so is
$\W_{\K_{\infty},p}/s\W_{\K_{\infty},p}$. This is a contradiction.
\end{proof}
This leads to the following obvious corollary.

\begin{corollary}\label{c:items}
Under the condition that no place in $\cS$ splits completely over
$\K_{\Xi}/k$, every object involved in the exact sequences {\em
(\ref{e:snake1})} and {\em (\ref{e:snake2})} is a finitely generated
torsion $\Lambda_{\Upsilon/\Xi}$-module.
\end{corollary}

\begin{corollary}\label{c:inject}
Under the condition that no place in $\cS$ splits completely over
$\K_{\Xi}/k$, the map
$\mathcal{T}_{\K_{\infty}}^{\Xi}\longrightarrow
\image(\phi)/s\image(\phi)$ in the exact sequence {\em
(\ref{e:snake2})} is injective.
\end{corollary}
\begin{proof}
Since $s$ is relatively prime to every $\zeta_i$, the group
$(\bigoplus_{j=1}^J \Lambda_{\Upsilon}/(\zeta_j^{m_j}))^{\Xi}$ is in
fact trivial.
\end{proof}

From this corollary and the exact sequence (\ref{e:snake2}), we get
\begin{equation}\label{e:character2}
\chi_{{\Upsilon/\Xi}}( \image(\phi)/s\image(\phi))\cdot
\chi_{{\Upsilon/\Xi}}( \mathcal{T}_{\K_{\infty}}/s
\mathcal{T}_{\K_{\infty}})
=\chi_{{\Upsilon/\Xi}}(\mathcal{T}_{\K_{\infty}}^{\Xi})\cdot
\prod_{j=1}^J({\bar \zeta}_j^{m_j}).
\end{equation}
Also, from the injection (\ref{e:nsn}) and the exact sequence
(\ref{e:snake1}), we get
\begin{equation}\label{e:character1}
\chi_{{\Upsilon/\Xi}}(\W_{\K_{\infty},p}/s\W_{\K_{\infty},p})=
\chi_{{\Upsilon/\Xi}}(\image(\phi)/s\image(\phi))\cdot
\chi_{{\Upsilon/\Xi}}(\mathcal{N}_{\K_{\infty}}/s
\mathcal{N}_{\K_{\infty}}).
\end{equation}
We use the
following lemma to make further simplification of the above relations.

\begin{lemma}\label{l:reduction}
We have
$$\chi_{{\Upsilon/\Xi}}( \mathcal{T}_{\K_{\infty}}/s \mathcal{T}_{\K_{\infty}})
=\chi_{{\Upsilon/\Xi}}(\mathcal{T}_{\K_{\infty}}^{\Xi}),$$ and
$$\chi_{{\Upsilon/\Xi}}(\mathcal{N}_{\K_{\infty}}/s \mathcal{N}_{\K_{\infty}})
=\chi_{{\Upsilon/\Xi}}(\mathcal{N}_{\K_{\infty}}^{\Xi}).$$
\end{lemma}
\begin{proof}
We have the obvious exact sequence
$$0\longrightarrow \mathcal{T}_{\K_{\infty}}^{\Xi}\longrightarrow
\mathcal{T}_{\K_{\infty}}\stackrel{s}{\longrightarrow}
\mathcal{T}_{\K_{\infty}}\longrightarrow \mathcal{T}_{\K_{\infty}}/s
\mathcal{T}_{\K_{\infty}} \longrightarrow 0,$$ where each term
is actually a finitely generated torsion
$\Lambda_{\Upsilon/\Xi}$-module. The first equality is proved by using the
fact that the characteristic ideals are multiplicative. The second
equality is proved by a similar argument.
\end{proof}

Finally we  prove the following key lemma.

\begin{lemma}\label{l:finally}
Suppose that $\K_{\infty}/k$ is a
$\Z_p^e$-extension with $\Gal(\K_{\infty}/k)=\Upsilon$ and $\Xi$  is a
 rank one $\Z_p$-submodule of $\Upsilon$ with $\Upsilon/\Xi\simeq
\Z_p^{e-1}$. Let $\K_{\Xi}$ be the fixed field of $\Xi$.
We define $\W_{\K_{\infty},p}$  and  $\W_{\K_{\Xi},p}$
as in {\em (\ref{e:wkinfp})}  and $\fG_v$ as in {\em (\ref{e:wvim})}.
Let $\cS$ be a set of places of $k$ so that
$\K_{\infty}/k$ is unramified outside of $\cS$ and
no place in $\cS$ splits
completely in $\K_{\Xi}/k$.
Let $\cS_1=\{v\in \cS \;\mid\;\Upsilon_v^0 \cap\Xi\not=\{0\}\}
$ and
$\cS_2=\{v\in \cS_1\;\mid\; |   \Upsilon_v/\Upsilon_v^0|<\infty\}$.
Let ${\bar \zeta}$ denotes the image of $\zeta$ under the projection
$\Lambda_{\Upsilon}\longrightarrow\Lambda_{\Upsilon/\Xi}.$
Then
$$
\chi_{{\Upsilon/\Xi}}(\W_{\K_{\Xi},p}) \cdot
\prod_{v\in\cS_1\setminus\cS_2}\chi_{{\Upsilon/\Xi}} (\fG_v)=
{\overline{\chi_{{\Upsilon}}(\W_{\K_{\infty},p})}} \cdot
\chi_{{\Upsilon/\Xi}}(\Z_p).
$$
Here the module $\Z_p$ is endowed with the trivial action of
$\Upsilon$.
\end{lemma}

\begin{proof}
Lemma \ref{l:reduction},  Corollary \ref{c:invxi}, equations (\ref{e:character2}) and
(\ref{e:character1}) imply that
$$\chi_{{\Upsilon/\Xi}}(\W_{\K_{\infty},p}/s\W_{\K_{\infty},p})=
\prod_{v\in\cS_2}\chi_{\Upsilon/\Xi} (\fG_v) \cdot
\prod_{j=1}^J({\bar \zeta}_j^{m_j}).
$$ From Lemma \ref{l:summerize}, we obtain
$$
\chi_{{\Upsilon/\Xi}}(\W_{\K_{\infty},p}/s\W_{\K_{\infty},p}) \cdot
\chi_{{\Upsilon/\Xi}}(\Z_p)=
\chi_{{\Upsilon/\Xi}}(\W_{\K_{\Xi},p})\cdot\prod_{v\in\cS_1}\chi_{\Upsilon/\Xi}
(\fG_v),$$
and this gives the lemma.
\end{proof}
\end{subsection}

\end{section}

\begin{section}{Proof of the Main Theorem}\label{s:end}

 The first step of the proof is to apply Lemma \ref{l:indp}. It
allows us to find an independent extension $(\L_{\infty}/k,{\tilde
S})$ of the given pair $(k_{\infty}/k,S)$ with $|\tilde S|$ much
larger than $|S|$. Let $k'_{\infty}/k$ be a $\Z_p$-extension
satisfying Lemma \ref{l:indp} (3). Since $k'_{\infty}/k$ is ramified
at every place in $\tilde S$ while the extension $k_{\infty}/k$ is
unramified at each $v\in {\tilde S}\setminus S$, we must have
$k'_{\infty}\nsubseteq k_{\infty}$. Therefore, the Galois group
$\Gal(k_{\infty}k'_{\infty}/k)$ is isomorphic to $\Z_p^{d+1}$.

Denote $\Gal(\L_{\infty}/k)=\Delta\simeq \Z_p^c$. By Lemma \ref{l:indp}, the
Stickelberger element $\theta_{\L_{\infty}/k,{\tilde S},T}$ is an
irreducible element in $\Lambda_{\Delta}$ and hence by Corollary
\ref{c:tate} and Corollary \ref{c:charw}, we have
\begin{equation}\label{e:m}
\chi_{{\Delta}}(\W_{\L_{\infty},p})=(\theta_{\L_{\infty}/k,{\tilde
S},T}^m),\; \text{for some}\; m\in\Z_+.
\end{equation}

To apply the results obtained in the previous sections, we
choose a chain of ascending Galois groups
\begin{equation}\label{e:seq2}
\Xi_0\subset\Xi_1\subset \dots \subset
\Xi_{c-1}=\Gal(\L_{\infty}/k'_{\infty})
\end{equation}
such that $\Xi_0=0$, $\Xi_{i}/\Xi_{i-1}\simeq\Z_p$ for $i>0$.

Let $\L_{\Xi_i}$, $i=0,..,c-1$, be
the fixed fields of $\Xi_i$ acting on $\L_{\infty}$.
Then we have
$$k \subset k'_{\infty}=\L_{\Xi_{c-1}}\subset\dots \subset \L_{\Xi_{1}}
\subset \L_{\Xi_{0}}=\L_{\infty}.$$
We shall also make the choice so that
\begin{equation}\label{e:both}
\L_{\Xi_{c-d-1}}=k_{\infty}k'_{\infty}.
\end{equation}

Now consider the sequence (\ref{e:seq2}) and for $i=0,...,c-2$ set
$\cS=\tilde S$, $\Upsilon=\Delta/\Xi_i$, $\Xi=\Xi_{i+1}/\Xi_{i}$,
$\K_{\infty}=\L_{\Xi_i}$ and $\K_{\Xi}=\L_{\Xi_{i+1}}$. We first note that the field
$\K_{\Xi}$ always contains the field $k'_{\infty}$, and since the
extension $k'_{\infty}/k$ is monomial it is ramified at every place in $\cS$.
In particular, no
place in $\cS$ splits completely in  $\K_{\Xi}/k$ and we can apply the key
Lemma \ref{l:finally} to these cases.  Also, if a place $v\in\cS_1\setminus\cS_2$, then
the corresponding residue extension for $\K_{\Xi}/k$ is of
infinite degree, hence the $\Z_p$-rank of decomposition subgroup
$(\Upsilon/\Xi)_v$ of $\Upsilon/\Xi$ is at least $2$. But
for $v\in\cS_1$ we have
$\fG_v=\Lambda_{\Upsilon}/\Im_v\simeq \Lambda_{\Upsilon/\Xi}/ {\bar \Im}_v$
where ${\bar \Im}_v$ is the ideal of $\Lambda_{\Upsilon/\Xi}$
generated by the set
$\{\sigma-1\; | \; \sigma\in (\Upsilon/\Xi)_v\}$.
Therefore, for $v\in\cS_1\setminus\cS_2$, the module $\fG_v$ is in
fact pseudo-null over $\Lambda_{\Upsilon/\Xi}$.

We also note that in the case where $i\leq c-3$, the $\Z_p$-rank of
$\Upsilon/\Xi$ is at least $2$. Consequently, the module $\Z_p$ is
pseudo-null over $\Lambda_{\Upsilon/\Xi}$ which is the ring of
formal power series in at least two variables. Therefore its
characteristic ideal is $\Lambda_{\Upsilon/\Xi}$. Lemma \ref{l:finally} implies that
if we choose for each $i$ a generator $\zeta_i$ for the
characteristic ideal $\chi_{\Delta/\Xi_i}(\W_{\L_{\Xi_{i}},p}),$
then for $i=0,...,c-3$, the generator $\zeta_{i+1}$ can be chosen
as ${\bar \zeta_i}$ which is the image of $\zeta_i$ under the
projection
$\Lambda_{\Delta/\Xi_i}\longrightarrow \Lambda_{\Delta/\Xi_{i+1}}$.
Here we should remind the readers that according to the functoriality
(\ref{e:functo}), this ring homomorphism actually sends the
Stickelberger element $\theta_{\L_{\Xi_i}/k,{\tilde S},T}$ to the
corresponding $\theta_{\L_{\Xi_{i+1}}/k,{\tilde S},T}$. This simple
fact turns out very useful, since the equation (\ref{e:m}) says that
$\zeta_0$ can be chosen as $\theta_{\L_{\Xi_0}/k,{\tilde S},T}^m$.
And from this we can deduce step by step that for $j=1,...,c-2$,
\begin{equation}\label{e:j}
\chi_{\Delta/\Xi_{j}}(\W_{\L_{\Xi_{j}},p})=(\zeta_{j})=
(\theta_{\L_{\Xi_{j}}/k,{\tilde S},T}^m).
\end{equation}
In the case where $i=c-2$ the situation is a little different from
the above. This time the Galois group
$\Upsilon/\Xi=\Delta/\Xi_{c-1} =\Gal(k'_{\infty}/k)$ is of rank one
over $\Z_p$. Then Lemma \ref{l:finally} says that
$$\chi_{{\Gal(k'_{\infty}/k)}}(\W_{k'_{\infty},p})=(\zeta_{c-1})=(
{\bar \zeta}_{c-2})\cdot \chi_{{\Gal(k'_{\infty}/k)}}(\Z_p).$$
Since $\Gal(k'_{\infty}/k)$ acts trivially on $\Z_p$, if $\sigma'$
is a topological generator of this Galois group then the ideal
$\chi_{{\Gal(k'_{\infty}/k)}}(\Z_p)=(\sigma'-1)$.
>From this and Corollary \ref{c:charw}, we get
$\chi_{{\Gal(k'_{\infty}/k)}}(\Cl_{k'_{\infty},p})=(
{\bar \zeta}_{c-2})$.  Then equation (\ref{e:j}), for $j=c-2$, and
the functoriality (\ref{e:functo}) allow us to conclude that
\begin{equation}\label{e:c-1}
\chi_{{\Gal(k'_{\infty}/k)}}(\Cl_{k'_{\infty},p})=
(\theta_{k'_{\infty}/k,{\tilde S},T}^m).
\end{equation}

On the other hand, the condition  we set at the beginning that
$k'_{\infty}/k$ is monomial implies that
$$\chi_{{\Gal(k'_{\infty}/k)}}(\Cl_{k'_{\infty},p})=
(\theta_{k'_{\infty}/k,{\tilde S},T})=(\sigma'-1)^{\tilde r}$$ where
${\tilde r}$ is  $|\tilde S|  -1$.
Comparing this with the equation (\ref{e:c-1}) and taking into account the
fact that the element $\sigma'-1$ is irreducible in
$\Lambda_{\Gal(k'_{\infty}/k)}$, we deduce that
$m=1$.

To complete the proof, we apply a similar argument
by setting
$\cS=\tilde S$, $\Upsilon=\Gal(k_{\infty}k'_{\infty}/k)$,
$\Xi=\Gal(k_{\infty}k'_{\infty}/k_{\infty})$,
$\K_{\infty}=k_{\infty}k'_{\infty}$ and $\K_{\Xi}=k_{\infty}$.
>From the equation
(\ref{e:both}) and (\ref{e:j}), with $j=c-d-1$, we already have
\begin{equation}\label{e:k'k}
\chi_{{\Gal(k'_{\infty}k_{\infty}/k)}}(\W_{k'_{\infty}k_{\infty},p})=(\theta_{
k'_{\infty}k_{\infty}/k,{\tilde S},T}).
\end{equation}
It is known (Corollary \ref{c:clks}) that no place in $\cS$ splits completely
in $\K_{\Xi}/k=k_{\infty}/k$.
Let $S_0\subset S$ be
the subset consisting of unramified places in $k_{\infty}/k$
and let $S_1=(S\cap \cS_1)\setminus (S_0\cup \cS_2)$. It is easy to see that
$\cS_1\setminus\cS_2=S_0  \sqcup S_1\sqcup ({\tilde S}\setminus S)$ and if $v\in S_1$ then
$\fG_v$ is pseudo-null over $\Lambda_{\Upsilon/\Xi}$ (the decomposition group $(\Upsilon/\Xi)_v$ is
at least of rank two over $\Z_p$).
In this situation, equations (\ref{e:functo}), (\ref{e:k'k}) and
Lemma \ref{l:finally} imply that
$$\chi_{{\Gamma}}(\W_{k_{\infty},p})\cdot
\prod_{v\in S_0 \sqcup({\tilde S}\setminus S)} (1-[v])
=(\theta_{k_{\infty}/k,{\tilde S},T})\cdot \chi_{{\Gamma}}(\Z_p),$$
where $[v]\in\Gamma$ is the Frobenius at $v$.
By
(\ref{e:compare}), this  reduces  to
$$\chi_{{\Gamma}}(\W_{k_{\infty},p})\cdot \prod_{v\in S_0} (1-[v])
=(\theta_{k_{\infty}/k,S,T})\cdot \chi_{{\Gamma}}(\Z_p).$$ If $d\geq
2$, then $\chi_{{\Gamma}}(\Z_p)=\Lambda_{\Gamma}$ and the proof is
completed. If $d=1$, then $\chi_{{\Gamma}}(\Z_p)=(\sigma-1)$, where $\sigma$ is a generator
of $\Gamma$, and we
can apply Corollary \ref{c:charw}.
This completes the proof of the Main Theorem.

\end{section}

\end{document}